\documentclass[a4paper,11pt]{amsart}

\usepackage[headings]{fullpage}

\calclayout

\makeindex

\usepackage{amsfonts,amsmath,amsthm,amscd,amssymb,latexsym,euscript, enumerate}
\usepackage{epsfig}
\usepackage{flafter}
\usepackage[all,cmtip,line]{xy}
\usepackage{array}
\usepackage[english]{babel}
\usepackage{overpic}
\usepackage{subfig}
\usepackage{multirow}
\usepackage{microtype}
\usepackage{wrapfig}
\usepackage{tabularray}
\usepackage{longtable}
\usepackage{supertabular}
\usepackage{caption}
\usepackage{mathrsfs}
\usepackage{tikz}
\usetikzlibrary{positioning}
\usepackage{tikz-cd}
\usepackage{float}
\allowdisplaybreaks

\usepackage{graphicx}
\usepackage{xcolor}
\usepackage{hyperref}
\hypersetup{
    colorlinks=true,    
    linkcolor=blue,          
    citecolor=blue,      
    filecolor=blue,      
    urlcolor=blue           
}

\newtheorem{theorem}{Theorem}[section]
\newtheorem{lemma}[theorem]{Lemma}
\newtheorem{proposition}[theorem]{Proposition}

\newtheorem*{theorem*}{Theorem}

\theoremstyle{plain}

\theoremstyle{definition} 
\newtheorem{definition}[theorem]{Definition}
\newtheorem{definition-lemma}[theorem]{Definition-Lemma}

\newtheorem{example}[theorem]{Example}

\newtheorem{remark}[theorem]{Remark}

\newtheorem{notation}[theorem]{Notation}
\numberwithin{equation}{section}

\newcommand{\C}{\mathbb{C}}
\newcommand{\R}{\mathbb{R}}
\newcommand{\Z}{\mathbb{Z}}

\newcommand{\Q}{\mathbb{Q}}

\def\P{\mathbb{P}}

\DeclareMathOperator{\ord}{ord}

\DeclareMathOperator{\vol}{vol}

\def\mult{\operatorname{mult}}
\def\Null{\operatorname{Null}}

\def\Aut{\operatorname{Aut}}

\def\Supp{\operatorname{Supp}}

\def\Sing{\operatorname{Sing}}

\def\lct{\operatorname{lct}}

\usepackage{mathtools}

\DeclarePairedDelimiterX{\norm}[1]{\lVert}{\rVert}{#1}

\title[K-stability of singular del Pezzo surfaces]
{K-stability of del Pezzo surfaces with a single quotient singularity}

\begin{document}

\author[I.-K~Kim]{In-Kyun Kim}
\author[D.-W.~Lee]{Dae-Won Lee}
\address[In-Kyun Kim]{June E Huh Center for Mathematical Challenges, Korea Institute for Advanced Study, 85 Hoegiro Dongdaemun-gu, Seoul 02455, Republic of Korea.}
\email{soulcraw@kias.re.kr}
\address[Dae-Won Lee]{Department of Mathematics, Ewha Womans University, 52 Ewhayeodae-gil, Seodaemun-gu, Seoul 03760, Republic of Korea}
\email{daewonlee@ewha.ac.kr}

\thanks{We would like to thank the referees for their careful reading of the paper and for their helpful comments. The authors are supported by the National Research Foundation of Korea (No. RS-2025-00513064). The second author is partially supported by Basic Science Research Program through the National Research Foundation of Korea (NRF) funded by the Ministry of Education (No. RS-2023-00237440 and 2021R1A6A1A10039823) and by Samsung Science and Technology Foundation under Project Number SSTF-BA2302-03.}

\date{\today}

\begin{abstract}
In this paper, we study the K-stability of del Pezzo surfaces with a
unique quotient singularity whose minimal resolution has exceptional
locus $E_1+E_2$, where $E_1^2=-n$, $E_2^2=-m$, and $E_1\cdot E_2=1$ for $n,m\geq2$.
\end{abstract}

\maketitle


\section{Introduction}\label{sect:intro}
The Yau--Tian--Donaldson conjecture asserts that a smooth Fano variety
\(X\) admits a K\"ahler--Einstein metric if and only if \(X\) is
K-polystable. This conjecture was proved in
\cite{CDS15a,CDS15b,CDS15c,Tian15a,Tian15b}. Consequently, every
K-stable smooth Fano variety admits a K\"ahler--Einstein metric.

A smooth del Pezzo surface \(S\) is K-polystable if and only if
\(S\cong\P^2\), \(S\cong\P^1\times\P^1\), or \(S\) is obtained by
blowing up \(\P^2\) at \(k\) points in general position, where
\(3\leq k\leq8\) \cite{Tian90,TY87}. However, for a given (possibly singular) Fano variety \(X\), determining whether \(X\) is K-(poly)stable is a delicate problem. The \(\delta\)-invariant provides valuative criteria for
K-semistability and uniform K-stability. More precisely, by \cite{BJ20}, the inequality \(\delta(X)>1\) is equivalent to uniform K-stability of \(X\), and, in particular, implies that \(X\) is K-stable. The Abban--Zhuang theory developed in \cite{AZ22} provides a
systematic method for estimating the \(\delta\)-invariant. See
Section \ref{sect:prelim} for details.

The \(\delta\)-invariants of smooth del Pezzo surfaces were computed in \cite{ACCFKMGSSV, CZ19, PW18}. Although the K-stability of some Du Val del Pezzo surfaces was determined in \cite{Che08,OSS16}, the \(\delta\)-invariants of all Du Val del Pezzo surfaces were recently computed in a series of papers \cite{Den1, Den2, Den4, Den3} by using the Abban--Zhuang theory. At the local level, \cite{Err} classified the local \(\delta\)-invariants of smooth del Pezzo surfaces of degree \(2\), showing, in particular, that irrational values occur precisely at points lying on exactly one \((-1)\)-curve. Furthermore, the K-stability of surfaces obtained by blowing up \(\P(1,1,n)\) at \(k\) general smooth points, where \(k\leq n+4\), was studied in \cite{KW25}. More precisely, let \(S_n^k\) denote a surface obtained
by blowing up \(k\) general smooth points of \(\P(1,1,n)\). For any
\(n+1\) general smooth points \(p_1,\dots,p_{n+1}\), there exists a
unique curve \(C\in\left|\mathcal O_{\P(1,1,n)}(n)\right|\) passing through them. Its strict transform \(\widetilde C\) is a \((-1)\)-curve. Contracting \(\widetilde C\) gives a birational morphism \(\pi\colon S_n^{n+1}\to\overline S_n^{n+1}\). We refer the reader to \cite{CP20} for details. Moreover, these surfaces are connected by blowing up general smooth  points as follows.
\begin{figure}[H]
  \centering
  \begin{tikzpicture}
  [-,auto,node distance=1.5cm, thick,main node/.style={circle,draw,font=\sffamily \Large\bfseries}]
    \node[text=black] (1) {$S_n^0\coloneqq \mathbb{P}(1,1,n)$};
    \node[text=black] (2) [right=0.5cm of 1] {$S_n^1$};
    \node[text=black] (3) [right of=2] {$\cdots$};
    \node[text=black] (4) [right of=3] {$S_n^{n+1}$};
    \node[text=black] (5) [right of=4] {$\cdots$};
    \node[text=black] (6) [right of=5] {$S_n^{n+4}$};
    \node[text=black] (7) [below of=4] {$\overline{S}_n^{n+1}$};
   
    \path[every node/.style={font=\sffamily}]
      (2) edge[->] (1)
      (3) edge[->] (2)
      (4) edge[->] (3)
      (5) edge[->] (4)
      (6) edge[->] (5)
      (4) edge[->] (7);  
  \end{tikzpicture}
  \caption{Del Pezzo surfaces with a single singularity of type \(\frac{1}{n}(1,1)\)}\label{fig:one sing}
\end{figure}
The Hilbert series of \(\overline S_n^{n+1}\) agrees with that of a
degree \(n+1\) hypersurface in \(\P(1,1,1,n)\). Moreover,
\(\overline S_n^{n+1}\) can be obtained by blowing up \(\P^2\) at
\(n+1\) points on a line and then contracting the strict transform of the line. See Example \ref{ex:Hilb1} and Proposition \ref{prop:one line}.

\begin{theorem}[{\cite[Theorem B]{KW25}}]\label{thm:one negative curve}
    Let \(n\geq 4\) be a positive integer, and \(S\) a del Pezzo surface which is isomorphic to one of the surfaces in Figure \ref{fig:one sing}. Then \(S\) is K-stable if and only if \(S\cong S_n^{n+4}\).
\end{theorem}

In this paper, we generalize Theorem \ref{thm:one negative curve} by allowing a prescribed subset of the blow-up centers to lie on a member of
\(\left|\mathcal O_{\P(1,1,n)}(1)\right|\).
\begin{definition}
    Let \(n,m\geq2\), and let \(\pi_n\colon \widetilde S_1\to\P(1,1,n)\) be the blow-up of \(m\) distinct smooth points on a member \(\ell_1\in\left|\mathcal O_{\P(1,1,n)}(1)\right|\). Let \(L_1\) be the strict transform of \(\ell_1\). Contracting \(L_1\) gives a birational morphism \(\varphi_n\colon \widetilde S_1\to S_{n,m}^0\).
\begin{figure}[H]
    \centering
    \begin{tikzpicture}
    [-,auto,node distance=2cm, thick,main node/.style={circle,draw,font=\sffamily \Large\bfseries}]
      \node[text=black] (1) {$\P(1,1,n)$};
      \node[text=black] (2) [above of=1, right of=1] {$\widetilde S_1$};
      \node[text=black] (3) [below of=2, right of=2] {$S_{n,m}^0$};
     
      \path[every node/.style={font=\sffamily}]
        (2) edge[->] node [above, left=0.5mm, pos=0.3] {\(\pi_n\)} (1)
        (2) edge[->] node [right, pos=0.5] {\(\varphi_n\)} (3);
    \end{tikzpicture}
    \caption{Construction of \(S_{n,m}^0\)}\label{fig:construction}
\end{figure}
    For \(k\geq0\), let \(S_{n,m}^k\) denote a surface obtained by blowing up \(k\) distinct general smooth points of \(S_{n,m}^0\).
\end{definition}

Throughout the paper, \(S_{n,m}^k\) denotes a member corresponding
to a point of a fixed nonempty Zariski-open subset of the relevant
configuration space. We choose this open subset so that all
incidence, irreducibility, transversality, and automorphism conditions
used below hold simultaneously. In the case \(S_{2,2}^4\), we also
require that the branch quartic of the anticanonical double cover
contains no line component.

We note that the surface \(S_{n,m}^n\) can also be obtained by blowing up points on the surface \(\P(1,1,m)\). Let \(\ell_2\in |\mathcal{O}_{\P(1,1,m)}(1)|\) be a curve, and let \(\pi_m\colon \widetilde S_2\to \P(1,1,m)\) be a blow-up at \(n\) distinct smooth points on \(\ell_2\). By contracting the strict transform \(L_2\) of \(\ell_2\), we obtain the birational morphism \(\varphi_m\colon \widetilde S_2\to S_{m,n}^0\). Let \(f\colon S_{m,n}^m\to S_{m,n}^0\) be a blow-up at \(m\) distinct general smooth  points on \(S_{m,n}^0\). For a suitable choice of these points, the surface \(S_{m,n}^m\) is isomorphic to \(S_{n,m}^n\). For \(n,m\geq2\), the surfaces \(S_{n,m}^{k_1}\) and
\(S_{m,n}^{k_2}\) are connected by the following blow-up diagram.
When \(n,m\geq3\), this is explained in \cite[Corollary 3.5.2]{Kut}; the cases in which one of \(n\) and \(m\) equals \(2\) follow from the same elementary transformations.
\begin{figure}[H]
    \centering
    \begin{tikzpicture}
    [-,auto,node distance=1.5cm, thick,main node/.style={circle,draw,font=\sffamily \Large\bfseries}]
      \node[text=black] (1) {$S_{n,m}^0$};
      \node[text=black] (2) [right of=1] {$S_{n,m}^1$};
      \node[text=black] (4) [right of=2] {$\cdots$};
      \node[text=black] (5) [right of=4] {$S_{n,m}^{n-1}$};
      \node[text=black] (7) [below of=5, right=2.5cm] {$S_{n,m}^n\cong S_{m,n}^m$};
      \node[text=black] (8) [below=2.5cm of 1] {$S_{m,n}^0$};
      \node[text=black] (9) [right of=8] {$S_{m,n}^1$};
      \node[text=black] (10) [right of=9] {$\cdots$};
      \node[text=black] (11) [right of=10] {$S_{m,n}^{m-1}$};
      \node[text=black] (12) [right=0.5cm of 7] {$S_{n,m}^{n+1}$};
      \node[text=black] (13) [right=0.5cm of 12] {$S_{n,m}^{n+2}$};
     
      \path[every node/.style={font=\sffamily}]
        (2) edge[->] (1)
        (4) edge[->] (2)
        (5) edge[->] (4)
        (7) edge[->] (5)
        (9) edge[->] (8)
        (10) edge[->] (9)
        (11) edge[->] (10)
        (7) edge[->] (11)
        (12) edge[->] (7)
        (13) edge[->] (12);
    \end{tikzpicture}
    \caption{Del Pezzo surfaces with a single singularity of type \(\frac{1}{mn-1}(1,n)\)}\label{fig:two negative}
\end{figure}

\begin{remark}\label{rem:n+3}
    Assume that \(n\geq m\geq 2\). The surface \(S_{n,m}^{n+3}\) is still a singular del Pezzo surface with a singularity of type \(\frac{1}{mn-1}(1,n)\) if and only if \((n,m)\) is one of the following: \((2,2),(3,2)\) or \((4,2)\).
\end{remark}
The surfaces in Figure \ref{fig:one sing} have only one exceptional curve in their minimal resolution. However, the surfaces in Figure \ref{fig:two negative} have exactly two exceptional curves in the minimal resolution. See Remark \ref{rem:HJ fraction} for details.

Now, we present the main result of this paper.

\begin{theorem}\label{main theorem}
Let \(S\) be a del Pezzo surface which is isomorphic to one of the
surfaces appearing in Figure \ref{fig:two negative} or Remark \ref{rem:n+3}.
Then the following assertions hold:
\begin{enumerate}[(1)]
    \item The surface \(S\) is K-stable if and only if it is isomorphic
    to one of
    \[
        S_{2,2}^4,\quad S_{2,2}^5,\quad
        S_{3,2}^5,\quad S_{3,2}^6,\quad
        S_{4,2}^6,\quad S_{4,2}^7.
    \]
    \item The surface \(S\) is strictly K-semistable if and only if it
    is isomorphic to one of
    \[
        S_{2,2}^3,\quad S_{3,3}^5,\quad S_{5,2}^7.
    \]
    \item Every remaining surface is K-unstable. In particular,
    \(S_{4,3}^6\) is K-unstable.
\end{enumerate}
\end{theorem}

Theorem \ref{main theorem} is summarized in Table \ref{Table}. Throughout Table \ref{Table}, we assume that \(n\geq m\).

    \begin{longtable}{cclclcl}
    \caption{K-stability of \(S_{n,m}^k\)}\label{Table}\\
    \hline\\[-1.5 ex]
    \((n,m)\) && \(k\) && K-stability && Reference \\
    \hline\hline\\[-1.5 ex]
    \(n+m\geq 8\) && \(k\leq n+2\) && K-unstable && Theorem \ref{thm:solution set}\\
    \hline\\[-1.5 ex]
    \((2,2)\) && \(k\leq 2\) && K-unstable && \cite{Den1}\\
    \hline\\[-1.5 ex]
    \((2,2)\) && \(k=3\) && strictly K-semistable && \cite{OSS16,Den2}\\
    \hline\\[-1.5 ex]
    \((2,2)\) && \(k=4,5\) && K-stable && \cite{Che08,Den4,Den3}\\
    \hline\\[-1.5 ex]
    \((3,2)\) && \(k\leq 4\) && K-unstable && Theorem \ref{thm:solution set}\\
    \hline\\[-1.5 ex]
    \((3,2)\) && \(k=5,6\) && K-stable && Theorems \ref{thm:325 stable} and \ref{thm:326 stable}\\
    \hline\\[-1.5 ex]
    \((4,2)\) && \(k\leq 5\) && K-unstable && Theorem \ref{thm:solution set}\\
    \hline\\[-1.5 ex]
    \((4,2)\) && \(k=6,7\) && K-stable && Theorems \ref{thm:426 stable} and \ref{thm:427 stable}\\
    \hline\\[-1.5 ex]
    \((3,3)\) && \(k\leq 4\) && K-unstable && Theorem \ref{thm:solution set}\\
    \hline\\[-1.5 ex]
    \((3,3)\) && \(k=5\) && strictly K-semistable && Theorem \ref{thm:335-semistable}\\
    \hline\\[-1.5 ex]
    \((4,3)\) && \(k\leq 5\) && K-unstable && Theorem \ref{thm:solution set}\\
    \hline\\[-1.5 ex]
    \((4,3)\) && \(k=6\) && K-unstable && Theorem \ref{thm:436-unstable}\\
    \hline\\[-1.5 ex]
    \((5,2)\) && \(k\leq 6\) && K-unstable && Theorem \ref{thm:solution set}\\
    \hline\\[-1.5 ex]
    \((5,2)\) && \(k=7\) && strictly K-semistable && Theorem \ref{thm:527-semistable}\\
    \hline\\[-1.5 ex]
    \end{longtable}

We note that \(S_{n,m}^{n+1}\) can also be obtained by running a \(-K_S\)-minimal model program for some smooth rational surface \(S\) with big anticanonical divisor \(-K_S\). Such a surface \(S\) is a Mori dream space by \cite{TVAV09}. See Proposition \ref{prop:two lines} and Remark \ref{rem:Xu big} for details.

By \cite[Theorem 1.2]{Xu}, if a klt pair \((X,\Delta)\) with
\(-(K_X+\Delta)\) big admits an anticanonical model
\((Z,\Delta_Z)\), then the K-stability properties of
\((X,\Delta)\) and \((Z,\Delta_Z)\) coincide. Consequently, the
generalized K-stability of the minimal resolution with respect to its
big anticanonical divisor agrees with the K-stability of its
anticanonical model.

The rest of the paper is organized as follows. In Section \ref{sect:prelim}, we
recall the basic notation for quotient surface singularities, the
\(\alpha\)-, \(\beta\)-, and \(\delta\)-invariants, and the relevant
criteria for K-stability. We also review the theory of potential
pairs. In Section \ref{sect:main}, we prove the main theorem. We first
use Theorem \ref{thm:normalized volume} to show that
\(S_{n,m}^k\) is K-unstable for every \(k\leq n+2\) whenever
\(m+n\geq8\). When \(m+n\leq7\), we apply the Abban--Zhuang theory
to prove that \(S_{2,2}^4\), \(S_{2,2}^5\), \(S_{3,2}^5\),
\(S_{3,2}^6\), \(S_{4,2}^6\), and \(S_{4,2}^7\) are K-stable.
We also prove that \(S_{3,3}^5\) and \(S_{5,2}^7\) are strictly
K-semistable, whereas \(S_{4,3}^6\) is K-unstable.

\section{Preliminaries}\label{sect:prelim}

Throughout this paper, we work over the complex number field \(\C\). We first recall the basic definition of the cyclic quotient singularity and of the Hirzebruch--Jung fraction.

\subsection{Cyclic quotient singularities of surfaces}
\begin{definition}
Let \(r\geq2\), let \(a\) be an integer satisfying \(1\leq a<r\) and \(\gcd(a,r)=1\), and put \(\varepsilon\coloneqq\exp\left(\frac{2\pi i}{r}\right)\). Let \(\Z/r\Z\) act on \(\C^2\) by \((x,y)\longmapsto(\varepsilon x,\varepsilon^a y)\). The quotient singularity at the image of the origin is said to be
of type \(\frac1r(1,a)\).

Let \(S\) be an algebraic surface. A point \(q\in S\) is said to be
of type \(\frac1r(1,a)\) if its analytic germ is isomorphic to the
germ of the origin in \(\C^2/(\Z/r\Z)\).
\end{definition}
\begin{definition}
    Let \(r\) and \(a\) be coprime integers with \(r>a>0\). The \textit{Hirzebruch--Jung fraction} of \(\frac{r}{a}\) is the expression
    \begin{align*}
        \frac{r}{a}=a_1-\frac{1}{a_2-\frac{1}{a_3-\dots}}=[a_1,\dots,a_k].
    \end{align*}
\end{definition}

\begin{proposition}[{\cite[p.~89]{Kol07}}]\label{prop:HJ fraction}
Let \(S\) be an algebraic surface, and let \(p\in S\) be a cyclic quotient singularity of type
\(\frac1r(1,a)\), and write
\[
\frac ra=[a_1,\dots,a_k].
\]
Then the exceptional locus of the minimal resolution \(\pi\colon\widetilde S\to S\) is a chain \(E_1+\cdots+E_k\) satisfying \(E_i^2=-a_i\), \(E_i\cdot E_{i+1}=1\) and \(E_i\cdot E_j=0\) whenever \(\lvert i-j\rvert>1\).
\end{proposition}

\begin{remark}\label{rem:HJ fraction}
    The surfaces \(S_n^k\) in Figure \ref{fig:one sing} have a single quotient singular point of type \(\frac{1}{n}(1,1)\). Since the Hirzebruch--Jung fraction of \(\frac{n}{1}\) is \([n]\), the minimal resolution has a unique exceptional curve \(E\) with \(E^2=-n\) by Proposition \ref{prop:HJ fraction}.

    The surfaces \(S_{n,m}^k\) in Figure \ref{fig:two negative} have a
    unique quotient singular point of type \(\frac{1}{mn-1}(1,n)\). Since
\[
\frac{mn-1}{n}=[m,n],
\]
    the exceptional locus of the minimal resolution is a chain \(F_1+F_2\) such that
\[
F_1^2=-m,\qquad F_2^2=-n,\qquad F_1\cdot F_2=1.
\]
    When convenient, we relabel \(E_1\coloneqq F_2\) and \(E_2\coloneqq F_1\), so that \(E_1^2=-n\) and \(E_2^2=-m\).
\end{remark}

\subsection{Embedding models}\label{subsec:embedding}
In this subsection, we note that some of the surfaces in Figures
\ref{fig:one sing} and \ref{fig:two negative} admit realizations as
hypersurfaces or complete intersections in weighted projective spaces.
We also show that \(\overline S_n^{n+1}\) and \(S_{n,m}^{n+1}\)
can be constructed from \(\P^2\), as described in Propositions
\ref{prop:one line} and \ref{prop:two lines}.

We will use the following result from \cite[Subsection 4.4]{CT}.
\begin{proposition}[{\cite[Subsection 4.4]{CT}}]\label{prop:ci hyp}
    Let \(n\) be a positive integer. The surface \(S_{2n-1}^{2n+3}\) can be embedded in \(\P(1,1,n,n,2n-1)\) as a complete intersection of two degree \(2n\) hypersurfaces. The surface \(S_{2n}^{2n+4}\) can be embedded in \(\P(1,1,n,n+1)\) as a degree \(2n+2\) hypersurface. 
\end{proposition}

The following Hilbert series can be computed using \cite{RS03}.

\begin{example}\label{ex:Hilb1}
The Hilbert series of \(\overline S_n^{n+1}\) is
\[
\frac{1-t^{n+1}}{(1-t)^3(1-t^n)}.
\]
This agrees with the Hilbert series of a degree \(n+1\)
hypersurface in \(\P(1,1,1,n)\).
\end{example}

\begin{example}
The Hilbert series of \(S_{n,m}^0\) is
\[
\frac{1+t^n+t^{2n}+\cdots+t^{n(m-1)}}
{(1-t)^2(1-t^{nm-1})}
=
\frac{1-t^{nm}}
{(1-t)^2(1-t^n)(1-t^{nm-1})}.
\]
This agrees with the Hilbert series of a degree \(nm\)
hypersurface in \(\P(1,1,n,nm-1)\).
\end{example}

However, in general, we do not know whether \(S_{n,m}^k\) can be
embedded as a hypersurface or a complete intersection in a weighted
projective space.

\bigskip

It was proved in \cite{TVAV09} that every smooth rational surface with
big anticanonical divisor is a Mori dream space. In particular,
\cite[Theorem 4.3]{TVAV09} treats blow-ups of \(\P^2\) at points on
a reducible cubic, under suitable conditions on the numbers of points
lying on its irreducible components. Propositions \ref{prop:one line}
and \ref{prop:two lines} show that \(\overline S_n^{n+1}\) and
\(S_{n,m}^{n+1}\) arise as anticanonical models of such surfaces.

\begin{proposition}\label{prop:one line}
Let \(L\subset\P^2\) be a line, and let
\(p_1,\dots,p_{n+1}\) be distinct points on \(L\). Let
\(\varphi\colon\widetilde S\to\P^2\) be their blow-up.
Contracting the strict transform \(\widetilde L\) of \(L\) gives a
birational morphism \(\pi\colon\widetilde S\to S\).
Then \(S\cong\overline S_n^{n+1}\).
\end{proposition}
\begin{proof}
Let \(B_n\subset\mathbb F_n\) be the negative section, and let
\(C_n\subset\mathbb F_n\) be a section with \(C_n^2=n\).
Blow up \(n+1\) distinct points
\(q_1,\dots,q_{n+1}\in C_n\), and let
\(\widetilde C_n\) and \(\widetilde F_i\) be the strict transforms of
\(C_n\) and of the fibre through \(q_i\), respectively. Then we have \(\widetilde C_n^2=\widetilde F_i^2=-1\).
Contracting \(\widetilde C_n\) gives the minimal resolution of
\(\overline S_n^{n+1}\).

The curves
\(\widetilde F_1,\dots,\widetilde F_{n+1}\) are pairwise disjoint.
Contracting them gives a smooth rational surface \(Z\) satisfying \(\rho(Z)=n+2-(n+1)=1\). Therefore, \(Z\cong\P^2\).

The curve \(B_n\) meets each \(\widetilde F_i\) once. Consequently, its
image \(B_Z\subset Z\) satisfies \(B_Z^2=-n+(n+1)=1\).
Thus, \(B_Z\) is a line in \(\P^2\). Moreover, the contraction points
of \(\widetilde F_1,\dots,\widetilde F_{n+1}\) lie on \(B_Z\).
Therefore, the minimal resolution is the blow-up of \(\P^2\) at
\(n+1\) points on a line, and contracting its strict transform gives
\(\overline S_n^{n+1}\).
\end{proof}

\begin{proposition}\label{prop:two lines}
Let \(L_1,L_2\subset\P^2\) be distinct lines. Let
\(p_1,\dots,p_{n+1}\) be distinct points on \(L_1\), and let
\(q_1,\dots,q_{m+1}\) be distinct points on \(L_2\). Assume that none
of these points is equal to \(L_1\cap L_2\). Let
\(\pi\colon S_1\to\P^2\) be the blow-up of these
\(n+m+2\) points, and let
\(\varphi\colon S_1\to S_2\) be the birational morphism
contracting the strict transforms of \(L_1\) and \(L_2\).
Then \(S_2\cong S_{n,m}^{n+1}\).
\end{proposition}
\begin{proof}
Let \(B_n\subset\mathbb F_n\) be the negative section, let
\(C_n\subset\mathbb F_n\) be a section with \(C_n^2=n\), and let
\(F\) be a fibre. Blow up \(m\) distinct points
\(r_1,\dots,r_m\in F\setminus(B_n\cup C_n)\) and \(n+1\) distinct
points \(s_1,\dots,s_{n+1}\in C_n\), all chosen away from
\(F\cap C_n\).

Let \(E_i\) be the exceptional curve over \(r_i\), let
\(\widetilde F_j\) be the strict transform of the fibre through \(s_j\),
and let \(\widetilde C_n\) and \(\widetilde F\) be the strict transforms
of \(C_n\) and \(F\), respectively. Then
\[
E_1,\dots,E_m,\quad
\widetilde F_1,\dots,\widetilde F_{n+1},\quad
\widetilde C_n
\]
are pairwise disjoint \((-1)\)-curves.

Contracting these curves gives a smooth rational surface \(Z\) with \(\rho(Z)=n+m+3-m-(n+1)-1=1\). Therefore, \(Z\cong\P^2\).

The image \(L_1'\subset Z\) of \(B_n\) satisfies \((L_1')^2=-n+(n+1)=1\).
The image \(L_2'\subset Z\) of \(\widetilde F\) satisfies \((L_2')^2=-m+m+1=1\). Hence, \(L_1'\) and \(L_2'\) are distinct lines in \(\P^2\).

The centers corresponding to
\(\widetilde F_1,\dots,\widetilde F_{n+1}\) lie on \(L_1'\), while the
centers corresponding to
\(E_1,\dots,E_m,\widetilde C_n\) lie on \(L_2'\). None of these centers
is equal to \(L_1'\cap L_2'\). Therefore, the surface before
contracting the two negative curves is precisely the blow-up described
in the statement, and the resulting surface is \(S_{n,m}^{n+1}\).
\end{proof}

\subsection{K-stability}

In this subsection, we recall the definitions and criteria for
K-stability that will be used below.

Let \(X\) be a \(\Q\)-Fano variety, that is, a normal projective variety with klt singularities such that \(-K_X\) is an ample \(\Q\)-Cartier divisor.

\begin{definition}
    A \textit{test configuration} \((\mathcal{X},\mathcal{L})\) of the pair \((X,-K_X)\) consists of 
\begin{enumerate}[(1)]
\item a normal variety \(\mathcal X\) equipped with a
\(\mathbb G_m\)-action;

\item a flat \(\mathbb G_m\)-equivariant morphism \(\pi\colon\mathcal X\to\P^1\), where \(\mathbb G_m\) acts on \(\P^1\) by \(\left(t,[x:y]\right)\longmapsto[tx:y]\);

\item a \(\mathbb G_m\)-linearized, \(\pi\)-ample
\(\Q\)-line bundle \(\mathcal L\) on \(\mathcal X\), together with a
\(\mathbb G_m\)-equivariant isomorphism
\[
\left(
\mathcal X\setminus\pi^{-1}([0:1]),
\mathcal L|_{\mathcal X\setminus\pi^{-1}([0:1])}
\right)
\cong
\left(
X\times\left(\P^1\setminus\{[0:1]\}\right),
\operatorname{pr}_1^{\ast}(-K_X)
\right).
\]
\end{enumerate}
\end{definition}

We write \(\mathcal X_0 \coloneqq \pi^{-1}([0:1])\) and \(\mathcal X_\infty \coloneqq \pi^{-1}([1:0])\). A test configuration \((\mathcal X,\mathcal L)\) is said to be of
\textit{product type} if there exists a one-parameter subgroup \(\lambda\colon\mathbb G_m\to\Aut(X)\) and a \(\mathbb G_m\)-equivariant isomorphism \(\left(
\mathcal X\setminus\mathcal X_{\infty},
\mathcal L|_{\mathcal X\setminus\mathcal X_{\infty}}
\right)
\cong
\left(
X\times\mathbb A^1,
\operatorname{pr}_1^{\ast}(-K_X)
\right)\),
where \(\mathbb G_m\) acts diagonally through \(\lambda\) on the first
factor and by multiplication on the second factor.

A product type test configuration is called \textit{trivial} if
\(\lambda\) is the trivial one-parameter subgroup.

The \textit{Donaldson--Futaki invariant} \(\mathrm{DF}(\mathcal{X};\mathcal{L})\) of the test configuration \((\mathcal{X},\mathcal{L})\) is defined as
\begin{align*}
    \mathrm{DF}(\mathcal{X};\mathcal{L})\coloneqq \frac{1}{(-K_X)^n}\left(\mathcal{L}^n\cdot K_{\mathcal{X}/\P^1}+\frac{n}{n+1}\mathcal{L}^{n+1}\right),
\end{align*}
where \(n=\dim X\).

\begin{definition}
    A Fano variety \(X\) is said to be \textit{K-semistable} if \(\mathrm{DF}(\mathcal{X};\mathcal{L})\geq 0\) for every test configuration \((\mathcal{X},\mathcal{L})\), and it is called \textit{K-stable} if \(\mathrm{DF}(\mathcal{X};\mathcal{L})> 0\) for every nontrivial test configuration \((\mathcal{X},\mathcal{L})\). We say a Fano variety \(X\) is \textit{K-polystable} if it is K-semistable and
    \begin{align*}
        \mathrm{DF}(\mathcal{X};\mathcal{L})=0 \Longleftrightarrow (\mathcal{X},\mathcal{L}) \text{ is of product type}.
    \end{align*}
    We say that \(X\) is \textit{strictly K-semistable} if \(X\) is K-semistable but not K-polystable.
\end{definition}

For a prime divisor \(E\) over \(X\), there exists a projective birational morphism \(\varphi\colon Y\to X\) such that \(E\) is a prime divisor on \(Y\). The \textit{center} \(C_X(E)\) of \(E\) is the image \(\varphi(E)\) of \(E\) on \(X\). The \textit{log discrepancy} of \(X\) along \(E\) is defined as 
\[
A_X(E)\coloneqq 1+\ord_E(K_Y-\varphi^{\ast}(K_X)).
\]
More generally, we define the log discrepancy for a log pair \((X,\Delta)\), i.e., \(\Delta\) is an effective \(\Q\)-divisor on \(X\) such that \(K_X+\Delta\) is \(\Q\)-Cartier. For a projective birational morphism \(\varphi\colon Y\to X\), write $K_{Y}+\Delta_{Y}=\varphi^{\ast}(K_{X}+\Delta)$ for some divisor \(\Delta_Y\) on \(Y\), and let $E$ be a prime divisor on $Y$. Define the log discrepancy \(A_{X,\Delta}(E)\) of \(E\) with respect to \((X,\Delta)\) as
\begin{align*}
    A_{X,\Delta}(E)\coloneqq 1+\ord_E(K_Y-\varphi^{\ast}(K_X+\Delta)).
\end{align*}

For a prime divisor \(E\) over \(X\), the \textit{pseudoeffective threshold} \(\tau_X(E)\) of \(E\) with respect to \(-K_X\) is defined as 
\begin{align*}
    \tau_X(E)\coloneqq \sup\{\tau\in \R_{\geq 0}\mid \vol(\varphi^\ast(-K_X)-\tau E)>0\}.
\end{align*}
The \textit{\(S\)-invariant} \(S_X(E)\) is defined as 
\[
S_X(E)\coloneqq \frac{1}{\vol(-K_X)}\int_0^{\tau_X(E)} \vol(\varphi^{\ast}(-K_X)-tE)\mathrm{d}t,
\]
and the \textit{\(\beta\)-invariant} \(\beta_X(E)\) is defined as 
\[
\beta_X(E)\coloneqq A_X(E)-S_X(E).
\] 
We will often omit \(X\) and simply write \(\tau(-)\), \(S(-)\) and \(\beta(-)\) if no confusion is likely to arise.

Next, we recall the definitions of the \(\alpha\)-, \(\beta\)-, and \(\delta\)-invariants, which give criteria for K-stability.

The \textit{\(\alpha\)-invariant} \(\alpha(X)\) is defined in \cite{Tian87} as 
\[
\alpha(X)\coloneqq \sup \left \{\lambda \in \Q \biggm\vert 
        \begin{array}{l} 
         \text{the log pair} ~ (X, \lambda D) ~ \text{is log canonical for every} \\
         \text{effective } \mathbb{Q}\text{-divisor } D \text{ on } X \text{ with } D\sim_{\Q} -K_X  
        \end{array}  \right \}.
\]

\begin{theorem}[{\cite{Tian87}}]
Let \(X\) be a smooth Fano variety of dimension \(n\). If \(\alpha(X)>\frac{n}{n+1}\), then \(X\) admits a K\"ahler--Einstein metric.
\end{theorem}

We have the following valuative criteria for K-stability.
\begin{theorem}[{\cite{BX19, Fuj19, Li17}}]
Let \(X\) be a \(\Q\)-Fano variety. Then \(X\) is K-semistable if and only if \(\beta_X(E)\geq 0\) for every prime divisor \(E\) over \(X\). In particular, if there exists a prime divisor \(E\) over \(X\) such that \(\beta_X(E)<0\), then \(X\) is K-unstable.
\end{theorem}

Let \(m\) be a positive integer such that \(-mK_X\) is Cartier, and \(\{s_1,\dots,s_{N_m}\}\) any basis of \(H^0(X,-mK_X)\). Let \(D_1,\dots,D_{N_m}\) be the corresponding divisors. Then the divisor 
\begin{align*}
    D=\frac{1}{mN_m}(D_1+\cdots+D_{N_m})\sim_{\Q} -K_X
\end{align*}
is called an \textit{\(m\)-basis type divisor} of \(-K_X\).

Let 
\[
\delta_m(X)\coloneqq \sup \left \{\lambda \in \Q \biggm\vert 
        \begin{array}{l} 
         \text{the log pair} ~ (X, \lambda D) ~ \text{is log canonical for every} \\
         \mathrm{effective} ~\mathbb{Q}\text{-divisor} ~D ~\text{on} ~X ~\text{with} ~D ~\text{of} ~m\text{-basis type}  
        \end{array}  \right \}.
\]
The \textit{\(\delta\)-invariant} \(\delta(X)\) was first defined as \(\delta(X)\coloneqq \limsup_{m\to \infty} \delta_m(X)\), and later it was shown that the \(\limsup\) is in fact a \(\lim\) \cite[Theorem 4.4]{BJ20}. By \cite[Theorem C]{BJ20}, we have
\begin{align*}
    \delta(X)=\inf_E \frac{A_X(E)}{S_X(E)},
\end{align*} 
where the infimum is taken over all prime divisors \(E\) over \(X\). The \textit{local \(\delta\)-invariant} \(\delta_p(X)\) of \(X\) at a point \(p\in X\) is defined as 
\begin{align*}
    \delta_p(X)\coloneqq \inf_{E} \frac{A_X(E)}{S_X(E)},
\end{align*}
where the infimum is taken over all prime divisors \(E\) over \(X\) such that \(p\in C_X(E)\). We note that \(\delta(X)=\inf_{p\in X} \delta_p(X)\). The definition of the \(\delta\)-invariant for a Fano variety can be generalized to a pair \((X,\Delta)\) such that \(-(K_X+\Delta)\) is big \cite{Xu}.

The following theorem explains why the \(\delta\)-invariant is often called the \textit{stability threshold}.

\begin{theorem}[{\cite[Theorem B]{BJ20}}]\label{thm:delta}
Let \(X\) be a \(\Q\)-Fano variety. Then the following hold:
\begin{enumerate}[(1)]
    \item \(X\) is K-semistable if and only if \(\delta(X)\geq 1\).
    \item \(X\) is uniformly K-stable if and only if \(\delta(X)> 1\).
\end{enumerate}
\end{theorem}

\begin{theorem}[{\cite[Theorem 1.6]{LXZ22}}]\label{thm:K-stable-uniform}
Let \(X\) be a \(\Q\)-Fano variety. Then \(X\) is K-stable if and only if
\(X\) is uniformly K-stable.
\end{theorem}

If the automorphism group is finite, then K-stability and K-polystability are equivalent conditions.

\begin{theorem}[{\cite[Corollary 1.5]{ACCFKMGSSV}}]\label{thm:finite aut}
    If \(\Aut(X)\) is finite, then \(X\) is K-stable if and only if it is K-polystable.
\end{theorem}

We have the implications
\[
    X \text{ is K-stable}\Longrightarrow X \text{ is K-polystable} \Longrightarrow X \text{ is K-semistable}.
\]
By Theorems \ref{thm:delta} and \ref{thm:K-stable-uniform}, if
\(\delta(X)=1\), then \(X\) is K-semistable and is not K-stable. If,
in addition, \(\Aut(X)\) is finite, then Theorem \ref{thm:finite aut}
implies that \(X\) is not K-polystable. Consequently, \(X\) is strictly
K-semistable.

By \cite{BJ20}, the \(\alpha\)-invariant also admits the valuative expression
\begin{align*}
    \alpha(X)=\inf_{E} \frac{A_X(E)}{\tau_X(E)},
\end{align*} 
where the infimum is taken over all prime divisors \(E\) over \(X\). For a point \(p\in X\), we also use the local \(\alpha\)-invariant 
\[
\alpha_p(X)
\coloneqq \inf_E \frac{A_X(E)}{\tau_X(E)},
\] 
where the infimum is taken over all prime divisors \(E\) over \(X\) such that \(p\in C_X(E)\). 

We have the following inequalities, which give lower and upper bounds for \(\delta(X)\) in terms of \(\alpha(X)\).

\begin{theorem}[{\cite[Theorem A]{BJ20}}]\label{thm:alpha delta}
Let \(X\) be a \(\Q\)-Fano variety of dimension \(n\). Then we have
\[
\frac{n+1}{n}\alpha(X)
\leq
\delta(X)
\leq
(n+1)\alpha(X).
\]
Moreover, the same argument, after restricting the infima to prime
divisors \(E\) satisfying \(p\in C_X(E)\), gives
\[
\frac{n+1}{n}\alpha_p(X)
\leq
\delta_p(X)
\leq
(n+1)\alpha_p(X)
\]
for every closed point \(p\in X\).
\end{theorem}

The following theorem plays a crucial role in proving Theorem \ref{thm:solution set}.

\begin{theorem}[{\cite[Theorem 3]{Liu}}]\label{thm:normalized volume}
    Let \(G\) be a finite group, \(X\) a \(\Q\)-Fano variety of dimension \(n\) and \(p\in X\) a quotient singularity with local analytic model \(\C^n/G\). If \(X\) is K-semistable, then we have 
    \[
    (-K_X)^n\leq \frac{(n+1)^n}{|G|}.
    \]
\end{theorem}

\subsection{Abban--Zhuang theory for surfaces}

Due to the development of the Abban--Zhuang theory \cite{AZ22}, determining the K-stability has become more feasible for various classes of varieties. For a variety \(X\) and a point \(p\in X\), the Abban--Zhuang theory gives a lower bound for the local \(\delta\)-invariant \(\delta_p(X)\). In this subsection, we recall how to estimate the \(\delta\)-invariant for surfaces by using Abban--Zhuang theory. 

Let \(S\) be a del Pezzo surface with at worst klt singularities and \(E\) a prime divisor over \(S\), i.e., there exists a projective birational morphism \(\varphi\colon \widetilde S\to S\) such that \(E\) is a prime divisor on \(\widetilde S\). The \(S(W_{\bullet,\bullet}^E;q)\)-invariant for a point \(q\in E\) can be calculated as follows.

\begin{theorem}[{\cite[Corollary 1.109]{ACCFKMGSSV}}]
    Let \(S\) be a del Pezzo surface with at worst klt singularities and \(\varphi\colon \widetilde S\to S\) a projective birational morphism such that \(E\) is a prime divisor on \(\widetilde S\). Let \(\varphi^{\ast}(-K_S)-tE=P(t)+N(t)\) be the Zariski decomposition. For each point \(q\in E\), we let
    \begin{equation}\label{function h(t)}
    \begin{split}
        h(t)&\coloneqq (P(t)\cdot E)\cdot \ord_q(N(t)\vert_E)+\int_0^{\infty} \vol_E(P(t)\vert_E-uq)\mathrm{d}u\\
        &=(P(t)\cdot E)(N(t)\cdot E)_{q}+\frac{1}{2}(P(t)\cdot E)^2,
    \end{split}
    \end{equation}
    where \((N(t)\cdot E)_{q}\) is the local intersection number of \(N(t)\) and \(E\) at the point \(q\).
    
    Then we have
    \[
    S(W_{\bullet,\bullet}^E;q)=\frac{2}{(-K_S)^2}\int_0^{\tau(E)} h(t)\mathrm{d}t,
    \]
    where \(\tau(E)\) is the pseudoeffective threshold of \(E\) with respect to \(-K_S\).
\end{theorem}

\begin{theorem}[{\cite[Theorem 3.2]{AZ22}}]\label{thm:AZ-primitive}
Let \(X\) be a klt Fano surface, let \(p\in X\), and let \(F\) be a
primitive divisor over \(X\) with center \(p\). Let
\(\pi_F\colon Y_F\to X\) be the associated prime blow-up. Assume that
\(F\) is Cartier on \(Y_F\) or is of plt type. Let
\(W_{\bullet,\bullet}^F\) be the refinement of the complete
anticanonical graded linear series by \(F\), and let \(\Phi_F\) be the
different determined by
\[
    (K_{Y_F}+F)|_F=K_F+\Phi_F.
\]
Then
\[
    \delta_p(X)\geq
    \min\left\{
        \frac{A_X(F)}{S_X(F)},
        \inf_{q\in F}
        \frac{A_{F,\Phi_F}(q)}
        {S(W_{\bullet,\bullet}^F;q)}
    \right\}.
\]
\end{theorem}

\begin{remark}\label{rem:AZ-higher-model}
Let \((X,\Delta)\) be a klt surface pair, and let
\(\mu\colon Y\to X\) be a projective birational morphism from a
smooth surface. Write \(K_Y+B_Y=\mu^{\ast}(K_X+\Delta)\).
Suppose that \(F\subset Y\) is a smooth prime divisor. Write \(B_Y=\left(1-A_{X,\Delta}(F)\right)F+B_Y'\), where \(F\not\subseteq\Supp(B_Y')\).
The different relevant to the Abban--Zhuang estimate on the crepant
model \((Y,B_Y)\) is defined by \((K_Y+F+B_Y')|_F=K_F+\Phi_F\).
This different need not agree with the different obtained after
contracting the other exceptional curves and passing to the associated
prime blow-up of \(F\) over \(X\).

Suppose that a prime divisor \(G\subset Y\), distinct from \(F\),
meets \(F\) transversely at a point \(q\), that \(q\) lies on no
other irreducible component of \(B_Y'\), and that the coefficient of
\(G\) in \(B_Y'\) is \(1-A_{X,\Delta}(G)\). Then \(A_{F,\Phi_F}(q)=A_{X,\Delta}(G)\).
In particular, if \(G\) is nonexceptional over \(X\) and is not
contained in the birational transform of \(\Delta\), then \(A_{F,\Phi_F}(q)=1\).
\end{remark}

\begin{lemma}\label{lem:AZ-crepant-model}
Let \((X,\Delta)\) be a projective klt surface pair, let \(M\) be a
big \(\Q\)-Cartier divisor on \(X\), and let \(p\in X\) be a closed
point. Let \(\mu\colon Y\to X\)
be a projective birational morphism from a smooth surface. Write \(K_Y+B_Y=\mu^{\ast}(K_X+\Delta)\), \(M_Y\coloneqq\mu^{\ast}M\),
and assume that \(B_Y\geq0\).

For a prime divisor \(V\) over \(X\), choose a projective birational
morphism \(f\colon Z\to X\) such that $V$ appears as a prime divisor on $Z$, and define
\[
S(M;V)
\coloneqq
\frac{1}{\vol(M)}
\int_0^{\infty}
\vol\bigl(f^{\ast}M-tV\bigr)\,\mathrm dt.
\]
This number is independent of the choice of \(f\). We define
\[
\delta_p(X,\Delta;M)
\coloneqq
\inf_{V}
\frac{A_{X,\Delta}(V)}{S(M;V)},
\]
where the infimum is taken over all prime divisors \(V\) over \(X\) such that \(p\in C_X(V)\). For prime divisors over $Y$, we use the analogous notation $S(M_Y;V)$.

For a closed point \(q\in Y\), put
\[
    \delta_q(Y,B_Y;M_Y)
    \coloneqq
    \inf_V
    \frac{A_{Y,B_Y}(V)}
    {S(M_Y;V)},
\]
where the infimum is taken over all prime divisors \(V\) over \(Y\)
such that \(q\in C_Y(V)\). Then
\[
    \delta_p(X,\Delta;M)
    =
    \inf_{q\in\mu^{-1}(p)}
    \delta_q(Y,B_Y;M_Y).
\]

Let \(F\subset Y\) be a smooth prime divisor. Write \(B_Y=\left(1-A_{X,\Delta}(F)\right)F+B_Y'\), where \(F\not\subseteq\Supp(B_Y')\),
and define \(\Phi_F\) by \((K_Y+F+B_Y')|_F=K_F+\Phi_F\).
Let \(W_{\bullet,\bullet}^{F}\) be the refinement of the complete
graded linear series of \(M_Y\) by \(F\). Then, for every closed point
\(q\in F\), we have
\[
    \delta_q(Y,B_Y;M_Y)
    \geq
    \min\left\{
        \frac{A_{X,\Delta}(F)}
        {S(M;F)},
        \frac{A_{F,\Phi_F}(q)}
        {S(W_{\bullet,\bullet}^{F};q)}
    \right\}.
\]
\end{lemma}

\begin{proof}
Let \(r\) be a sufficiently divisible positive integer such that
\(rM\) is Cartier. Since \(X\) is normal and
\(\mu_{\ast}\mathcal O_Y=\mathcal O_X\), the projection formula gives \(H^0(X,mrM) \cong H^0(Y,mrM_Y)\) for every \(m\geq0\). Under this identification, all divisorial
filtrations agree. Consequently, for every prime divisor \(V\) over
\(Y\), viewed also as a prime divisor over \(X\), we have \(A_{Y,B_Y}(V)=A_{X,\Delta}(V)\) and \(S(M_Y;V)=S(M;V)\). Moreover, \(\mu\bigl(C_Y(V)\bigr)=C_X(V)\). It follows that
\[
p\in C_X(V)
\quad\Longleftrightarrow\quad
C_Y(V)\cap\mu^{-1}(p)\neq\varnothing.
\]
Taking the corresponding infima gives
\[
\delta_p(X,\Delta;M)
=
\inf_{q\in\mu^{-1}(p)}
\delta_q(Y,B_Y;M_Y).
\]

Since \((X,\Delta)\) is klt and
\(K_Y+B_Y=\mu^{\ast}(K_X+\Delta)\), the crepant pair
\((Y,B_Y)\) is klt. The assumption \(B_Y\geq0\) therefore allows us to
apply the Abban--Zhuang estimate to this pair.

For the second assertion, apply
\cite[Theorem 3.2]{AZ22} to the pair \((Y,B_Y)\), the Veronese
complete graded linear series associated with \(rM_Y\), the Cartier
divisor \(F\), and the closed subset \(\{q\}\). Using the degree one homogeneity \(S(rM_Y;V)=rS(M_Y;V)\) and the analogous homogeneity for the restricted graded linear series, and then dividing by \(r\), we obtain
\[
\delta_q(Y,B_Y;M_Y)
\geq
\min\left\{
\frac{A_{Y,B_Y}(F)}{S(M_Y;F)},
\delta_q(F,\Phi_F;W_{\bullet,\bullet}^F)
\right\}.
\]
Since \(A_{Y,B_Y}(F)=A_{X,\Delta}(F)\) and \(S(M_Y;F)=S(M;F)\), it remains to compute the local threshold on \(F\). Since \(F\) is
a smooth curve, every divisorial valuation of \(\C(F)\) centered at
\(q\) is proportional to \(\ord_q\). Therefore,
\[
\delta_q(F,\Phi_F;W_{\bullet,\bullet}^F)
=
\frac{A_{F,\Phi_F}(q)}
{S(W_{\bullet,\bullet}^F;q)}.
\]
This proves the assertion.
\end{proof}

\subsection{Potential pairs}
Let \(Y\) be a smooth projective variety and \(D\) a big divisor on \(Y\). Let \(E\) be a prime divisor on \(Y\). We define the \textit{asymptotic divisorial valuation} \(\sigma_E(D)\) of \(D\) along \(E\) as \(\inf\{\mult_E(D')\mid 0\leq D'\sim_{\Q} D\}\). If \(D\) is only pseudoeffective, then we define \(\sigma_E(D)\coloneqq \lim_{\epsilon\to 0^{+}} \sigma_E(D+\epsilon A)\), where \(A\) is an ample divisor on \(Y\). It is well known that this definition is independent of the choice of an ample divisor. 

This definition naturally extends to normal projective varieties as follows. Let \(X\) be a normal projective variety and \(D\) a pseudoeffective divisor on \(X\). Let \(E\) be a prime divisor over \(X\). We define the asymptotic divisorial valuation \(\sigma_E(D)\coloneqq \sigma_E(f^{\ast}D)\), where \(f\colon Y\to X\) is a projective birational morphism from a smooth projective variety \(Y\) such that \(E\) is a prime divisor on \(Y\). This definition is independent of the choice of \(f\) as shown in  \cite[Theorem III.5.16]{Nak}.

For a normal projective variety \(X\) and a pseudoeffective divisor \(D\), we define the \textit{divisorial Zariski decomposition} as follows. We define the negative part \(N(D)\coloneqq \sum_{E} \sigma_E(D)E\), where the sum is taken over all prime divisors \(E\) on \(X\), and \(P(D)\coloneqq D-N(D)\) is defined as the positive part. We note that \(\sigma_E(D)>0\) for only finitely many prime divisors. The decomposition \(D=P(D)+N(D)\) is called the \textit{divisorial Zariski decomposition}. When \(\dim X\geq 3\), the positive part is not necessarily nef, but only movable, unlike the Zariski decomposition for surfaces. If the positive part is nef, then we simply call it the Zariski decomposition. See \cite{Nak} for more details.

\begin{definition}[{\cite[Definitions 3.1 and 3.2]{CP}}]
    Let \((X,\Delta)\) be a pair such that \(-(K_X+\Delta)\) is pseudoeffective. For a prime divisor \(E\) over \(X\), we define the \textit{potential discrepancy}  \(\overline{a}(E;X,\Delta)\) of \((X,\Delta)\) along \(E\) as
    \begin{align*}
        \overline{a}(E;X,\Delta)\coloneqq A_{X,\Delta}(E)-\sigma_E(-(K_X+\Delta)).
    \end{align*} 
    The pair \((X,\Delta)\) is called \emph{potentially klt}, or
\emph{pklt}, if
\(\inf_E\overline a(E;X,\Delta)>0\). It is called
\emph{potentially lc}, or \emph{plc}, if
\(\inf_E\overline a(E;X,\Delta)\geq0\), where in both cases the
infimum is taken over all prime divisors \(E\) over \(X\).
\end{definition}

The notion of potential pairs was introduced and studied in \cite{CP}. See also \cite{CJK, CJKL, CJL, Lee} for related results. One of the main properties of potential pairs is that one can bound the singularities of the pairs obtained by running a \(-(K_X+\Delta)\)-minimal model program (MMP for short). 

We briefly recall the notions of a \(-(K_X+\Delta)\)-minimal model program and an anticanonical model.
\begin{definition}
    Let \(X\) be a normal projective variety and \(D\) a pseudoeffective divisor on \(X\). A birational contraction \(\varphi\colon X\dashrightarrow Y\) is called \textit{\(D\)-negative} if \(\varphi_{\ast}D\) is \(\R\)-Cartier and there exists a common resolution \((p,q)\colon W\to X\times Y\) such that
    \begin{align*}
        p^{\ast}D=q^{\ast}\varphi_{\ast}D+E,
    \end{align*}
    where \(E\) is an effective \(q\)-exceptional divisor whose support \(\Supp(E)\) contains the support of all the strict transforms of the \(\varphi\)-exceptional divisors.
\end{definition}
As a special case of a \(D\)-negative contraction, we say \(\varphi\colon X\dashrightarrow Y\) is a \(D\)-minimal model if \(\varphi\) is a \(D\)-negative contraction and \(\varphi_{\ast}D\) is nef. The composition of \(D\)-negative contractions is called a \(D\)-minimal model program (\(D\)-MMP for short). When \(D=-(K_X+\Delta)\) is pseudoeffective, a \(D\)-minimal model and a \(D\)-MMP are called a \(-(K_X+\Delta)\)-minimal model and a \(-(K_X+\Delta)\)-MMP, respectively. 

\begin{definition}
    A rational map $f\colon X\dashrightarrow Y$ is called the \textit{anticanonical model} of $(X,\Delta)$ if $Y$ is a normal projective variety and there is an ample divisor $A$ on $Y$ such that if $(p,q)\colon W\rightarrow X \times Y$ is a resolution of the indeterminacy of $f$, then $q$ is a contraction morphism with $-p^\ast (K_X+\Delta)\sim_{\R} q^\ast A+E$, where $E\geq 0$ is contained in the fixed part of $\vert -p^\ast (K_X+\Delta)\vert_{\R}$. 
\end{definition}

\begin{theorem}[{\cite[Proposition 3.11]{CP}}]\label{thm:pklt klt}
    Let \((X,\Delta)\) be a pair such that \(-(K_X+\Delta)\) is pseudoeffective. Suppose that \(\varphi\colon (X,\Delta)\dashrightarrow (Y,\Delta_Y)\) is a \(-(K_X+\Delta)\)-minimal model. Then \((Y,\Delta_Y)\) is klt (resp. lc) if and only if \((X,\Delta)\) is pklt (resp. plc).
\end{theorem}

However, being a pklt pair does not guarantee that we can run an \(-(K_X+\Delta)\)-MMP in general. On the other hand, if \((X,\Delta)\) is a pklt pair such that \(-(K_X+\Delta)\) is big, then \(X\) is a variety of Fano type \cite[Theorem 5.1]{CP} and hence, we can run the \(-(K_X+\Delta)\)-MMP.

In the surface case, the positive part of the divisorial Zariski decomposition is nef. We will use the negative part of the Zariski decomposition to compute potential discrepancies.

\begin{proposition}[{\cite[Lemma 3.5, Corollary 5.2]{CP}}]\label{prop:pklt surface}
    Let \(S\) be a surface such that \(-K_S\) is big, and let \(-K_S=P+N\) be its Zariski decomposition. If \((S,N)\) is a klt pair, then \(S\) is of Fano type and \((S,0)\) is a pklt pair. In particular, by contracting all the curves in \(\Supp(N)\), we obtain the \(-K_S\)-minimal model.
\end{proposition}

For a klt pair \((X,\Delta)\) such that \(-(K_X+\Delta)\) is big, if one can run a \(-(K_X+\Delta)\)-MMP and obtain the anticanonical model \((Z,\Delta_Z)\) of \((X,\Delta)\), then the K-stability of the pairs \((X,\Delta)\) and \((Z,\Delta_Z)\) coincides, as shown by the following theorem.

\begin{theorem}[{\cite[Theorem 1.2]{Xu}}]\label{thm:Xu big}
    Let \((X,\Delta)\) be a klt pair such that \(-(K_X+\Delta)\) is big. Suppose that the anticanonical model \((Z,\Delta_Z)\) of \((X,\Delta)\) exists. Then \((X,\Delta)\) is K-semistable (resp. K-stable, uniformly K-stable) if and only if \((Z,\Delta_Z)\) is K-semistable (resp. K-stable, uniformly K-stable).
\end{theorem}

\begin{proposition}\label{prop:-K-mmp}
Let \(S_1\) and \(S_2\) be surfaces as in Proposition
\ref{prop:two lines}. Then \(S_2\) is the anticanonical model of
\(S_1\).

If \((n,m)\neq(2,2)\), then
\(\varphi\colon S_1\to S_2\) is \((-K_{S_1})\)-negative. If
\((n,m)=(2,2)\), then \(-K_{S_1}\) is nef and \(\varphi\) is the
crepant anticanonical contraction of \(\Null(-K_{S_1})\).
\end{proposition}
\begin{proof}
    On \(S_1\), let \(\widetilde L_n\) and \(\widetilde L_m\) be the strict transforms of \(L_1\) and \(L_2\), respectively. The anticanonical divisor \(-K_{S_1}\) can be written as \(-K_{S_1}=L+\widetilde L_n+\widetilde L_m\),
    where \(L\) is the pullback of the hyperplane class on \(\P^2\). The positive and negative parts of the Zariski decomposition of \(-K_{S_1}\) are as follows:
    \begin{align*}
        P=L+\frac{m+1}{mn-1}\widetilde L_n+\frac{n+1}{mn-1}\widetilde L_m, \qquad N=\frac{mn-m-2}{mn-1}\widetilde L_n+\frac{mn-n-2}{mn-1}\widetilde L_m.
    \end{align*}

If \((n,m)\neq(2,2)\), then we have \(-K_{S_1} = \varphi^{\ast}(-K_{S_2})+N\), and \(\varphi\) is \((-K_{S_1})\)-negative. Since the coefficients of \(N\) are smaller than \(1\), the pair \((S_1,N)\) is klt.
Proposition \ref{prop:pklt surface} implies that \((S_1,0)\) is pklt.

If \((n,m)=(2,2)\), then \(N=0\), and \(\Null(-K_{S_1}) = \widetilde L_n\cup\widetilde L_m\).
In this case, we have \(-K_{S_1}=\varphi^{\ast}(-K_{S_2})\). Thus, \(\varphi\) is crepant.

In either case, \(-K_{S_2}\) is ample. Therefore, \(S_2\) is the
anticanonical model of \(S_1\).
\end{proof}

\begin{remark}\label{rem:Xu big}
    By Proposition \ref{prop:two lines}, the surfaces \(S_{n,m}^{k}\) with \(k\geq n+1\) can be obtained by running an anticanonical MMP. More precisely, since they are singular del Pezzo surfaces, they are in fact the anticanonical models of surfaces obtained by blowing up \(\P^2\) in special configurations as in \cite{TVAV09}. Moreover, by Theorem \ref{thm:Xu big}, the generalized K-stability of the minimal resolution of \(S_{n,m}^k\), with respect to its big anticanonical divisor, agrees with the K-stability of \(S_{n,m}^k\).
\end{remark}

\section{Main results and proofs}\label{sect:main}

In this section, we prove Theorem \ref{main theorem}. Let us first show that when \(n+m\geq 8\), the surfaces \(S_{n,m}^{k_1}\) and \(S_{m,n}^{k_2}\) in Figure \ref{fig:two negative} are K-unstable by using Theorem \ref{thm:normalized volume}. Since we have the isomorphism \(S_{n,m}^n\cong S_{m,n}^m\), we need to consider the cases where \(k_1\leq n+2\) and \(k_2\leq m-1\).

\subsection{Case: \(n+m\geq 8\).}
In this section, we prove Theorem \ref{main theorem}. We first use
Theorem \ref{thm:normalized volume} to prove K-unstability when
\(n+m\geq8\).

\begin{lemma}\label{lem:inequalities}
    Assume that \(n\geq m\geq2\), and let \((n,m,k_1)\) and \((m,n,k_2)\) be triples occurring in Figure \ref{fig:two negative} or Remark \ref{rem:n+3}. The inequality
    \begin{align*}
        n+2-k_1+\frac{m+n+2}{mn-1}\leq \frac{9}{mn-1}
    \end{align*}
    holds if and only if the triple \((n,m,k_1)\) is contained in the following set:
    \begin{equation}\label{eq:exceptional-triples}
\mathscr A\coloneqq
\left\{
\begin{aligned}
&(2,2,3),(2,2,4),(2,2,5),(3,2,5),(3,2,6),\\
&(4,2,6),(4,2,7),(5,2,7),(3,3,5),(4,3,6)\\
\end{aligned}
\right\}.
\end{equation}
    Also, the inequality
    \begin{align*}
        m+2-k_2+\frac{m+n+2}{mn-1}\leq \frac{9}{mn-1}
    \end{align*}
    does not hold.
\end{lemma}
\begin{proof}
    For \(k_1\leq n+2\), the inequality can hold only if \(m+n\leq7\). Under \(n\geq m\geq2\), the possible pairs are \((2,2),(3,2),(4,2),(5,2),(3,3)\) and \((4,3)\). For these pairs, the smallest integer values of \(k_1\) satisfying the inequality are, respectively, \(3,5,6,7,5,6\). Together with the three \(k_1=n+3\) cases in Remark \ref{rem:n+3}, this gives precisely \(\mathscr A\).

    The second inequality is equivalent to
\[
m+2+\frac{m+n-7}{mn-1}\leq k_2.
\]
Since \(k_2\leq m-1\), it would follow that
\[
\frac{m+n-7}{mn-1}\leq-3.
\]
On the other hand, \(\frac{m+n-7}{mn-1}+1 = \frac{mn+m+n-8}{mn-1} \geq0\) for \(m,n\geq2\). Therefore, we obtain \(\frac{m+n-7}{mn-1}\geq-1\), which is a contradiction.
\end{proof}

\begin{theorem}\label{thm:solution set}
    Assume \(n\geq m\geq 2\). Let \(S_{n,m}^{k_1}\) and \(S_{m,n}^{k_2}\) be the surfaces as in Figure \ref{fig:two negative} and Remark \ref{rem:n+3}. The surface \(S_{m,n}^{k_2}\) is K-unstable for all \(n,m\geq 2\) and \(0\leq k_2\leq m-1\), and the surface \(S_{n,m}^{k_1}\) is K-unstable except when the triple \((n,m,k_1)\) belongs to the set \(\mathscr A\). In particular, if \(m+n\geq 8\), then the surface \(S_{n,m}^{k_1}\) is K-unstable for all \(k_1\leq n+2\).
\end{theorem}
\begin{proof}
    Let \(S_1\coloneqq S_{n,m}^{k_1}\) and \(S_2\coloneqq S_{m,n}^{k_2}\) for simplicity. We note that 
    \begin{align*}
        (-K_{S_1})^2=n+2-k_1+\frac{m+n+2}{mn-1}, \qquad (-K_{S_2})^2=m+2-k_2+\frac{m+n+2}{mn-1}.
    \end{align*}
    The singular points \(p\in S_1\) and \(q\in S_2\) are locally analytically isomorphic to \(\C^2/G_1\) and \(\C^2/G_2\) with \(|G_1|=|G_2|=mn-1\). Note that \(G_1\) acts on \(\C^2\) by \((x,y)\mapsto (\varepsilon x, \varepsilon^n y)\) and \(G_2\) acts on \(\C^2\) by \((x,y)\mapsto (\varepsilon x, \varepsilon^m y)\), where \(\varepsilon\coloneqq \exp\left(\frac{2\pi i}{mn-1}\right)\) is the \((mn-1)\)th primitive root of unity. If \(S_1\) is K-semistable, then Theorem \ref{thm:normalized volume} gives
\[
n+2-k_1+\frac{m+n+2}{mn-1}
\leq
\frac9{mn-1}.
\]
Similarly, if \(S_2\) is K-semistable, then
\[
m+2-k_2+\frac{m+n+2}{mn-1}
\leq
\frac9{mn-1}.
\]
Lemma \ref{lem:inequalities} now proves both assertions.
\end{proof}

\subsection{Case: \(m+n\leq 7\).}

Without loss of generality, we may assume that \(n\geq m\). By Theorem \ref{thm:solution set}, we only need to consider the surfaces \(S_{n,m}^k\) for which the triple \((n,m,k)\) is contained in the solution set \(\mathscr A\).

\begin{remark}\label{rem:22}
The \(\delta\)-invariants of the surfaces \(S_{2,2}^k\) were computed
in \cite{Den1,Den2,Den4,Den3}. Their unique singularity is of type
\(\frac13(1,2)\), hence is a Du Val singularity of type \(A_2\).
The surface \(S_{2,2}^k\) is K-unstable unless
\(k\in\{3,4,5\}\). More precisely, \(S_{2,2}^3\) is strictly
K-semistable, while \(S_{2,2}^4\) and \(S_{2,2}^5\) are K-stable.
K-stability for some of these surfaces had already been established
in \cite{Che08,OSS16}.
\end{remark}

The remaining cases are
\(S_{3,2}^5\), \(S_{3,2}^6\), \(S_{4,2}^6\), \(S_{4,2}^7\),
\(S_{5,2}^7\), \(S_{3,3}^5\), and \(S_{4,3}^6\).
We prove that the first four surfaces are K-stable, that
\(S_{3,3}^5\) and \(S_{5,2}^7\) are strictly K-semistable, and that
\(S_{4,3}^6\) is K-unstable.

\begin{lemma}\label{lem:base-22-smooth-alpha}
For every smooth point \(p\in S_{2,2}^4\) and every smooth point
\(q\in S_{2,2}^5\), we have
\[
\alpha_p(S_{2,2}^4)\geq\frac34
\qquad\text{and}\qquad
\alpha_q(S_{2,2}^5)\geq\frac57.
\]
In particular, both local \(\alpha\)-invariants are greater than
\(\frac7{10}\).
\end{lemma}

\begin{proof}
The surface \(S_{2,2}^4\) is a Gorenstein canonical del Pezzo surface of degree \(2\) with one \(A_2\)-singularity. Let \(\tau\colon S_{2,2}^4\to\P^2\) be the anticanonical double cover, and let \(p\in S_{2,2}^4\) be a smooth point.

Every member of \(|-K_{S_{2,2}^4}|\) is the inverse image of a line
in \(\P^2\). Let \(B\subset\P^2\) be the branch quartic of \(\tau\).
By the genericity convention imposed in the introduction, \(B\) has no
line component. If \(C=\tau^{-1}(\Lambda)\) is singular at \(p\),
where \(\Lambda\subset\P^2\) is a line, then the restriction of a
local equation of \(B\) to \(\Lambda\) has vanishing order
\(e\in\{2,3,4\}\) at \(\tau(p)\). Thus, in suitable analytic
coordinates at \(p\), the curve \(C\) has equation \(u^2=v^e\).
Consequently, \(\lct_p(S_{2,2}^4;C)\geq\frac34\). This is the degree-two local analysis used in
\cite[Section 4]{MG14}.

If \(C=\tau^{-1}(\Lambda)\) is reducible, say \(C=L_1+L_2\), then
each restriction
\(\tau|_{L_i}\colon L_i\to\Lambda\) is an isomorphism.
Consequently, each \(L_i\) is a smooth rational curve satisfying \((-K_{S_{2,2}^4})\cdot L_i=1\).

Suppose that there exist an effective \(\Q\)-divisor \(D\sim_{\Q}-K_{S_{2,2}^4}\) and a rational number \(\lambda<\frac34\) such that
\((S_{2,2}^4,\lambda D)\) is not log canonical at \(p\). Since \(p\) is smooth and
\((S_{2,2}^4,\lambda D)\) is not log canonical at \(p\), we have \(\mult_p(D)>\frac1\lambda\).
Put
\[
\mathcal C_p
\coloneqq
\{C\in|-K_{S_{2,2}^4}|\mid p\in C\}.
\]

Suppose first that some \(C\in\mathcal C_p\) is singular at \(p\).
Since \((S_{2,2}^4,\lambda C)\) is log canonical at \(p\), the convexity
argument in \cite[Lemma 2.4]{MG14} allows us to assume that an
irreducible component of \(C\) passing through \(p\) is not contained
in \(\Supp(D)\).

If \(C\) is irreducible, then
\[
2=C\cdot D
\geq
\mult_p(C)\mult_p(D)
>
\frac2{\lambda},
\]
which is impossible.

Suppose that \(C=L_1+L_2\) is reducible and that
\(L_1\not\subseteq\Supp(D)\), where \(p\in L_1\). Since \((-K_{S_{2,2}^4})\cdot L_1=1\),
we have
\[
1=L_1\cdot D
\geq
(L_1\cdot D)_p
\geq
\mult_p(D)
>
\frac1{\lambda},
\]
which is again impossible.

Therefore, every member of \(\mathcal C_p\) is smooth at \(p\).
Let \(\sigma\colon\widehat S_{2,2}^4\to S_{2,2}^4\) be the blow-up of \(p\), and let \(F\) be its exceptional curve.
Write \(\sigma^{\ast}D=D'+rF\), where \(r\coloneqq\mult_p(D)\).
The pair \(\left(\widehat S_{2,2}^4,\lambda D'+(\lambda r-1)F\right)\) is not log canonical at some point \(z\in F\).

Choosing a general member of \(\mathcal C_p\) that is not contained
in \(\Supp(D)\), we obtain \(r\leq2\).
By the preceding multiplicity inequality, \(\lambda r>1\). Applying
\cite[Lemma 2.5(i)]{MG14} to \(\widehat S_{2,2}^4\), we obtain
\[
\mult_z(D')+r>\frac2{\lambda}.
\]

Consider the differential map \(H^0\bigl(S_{2,2}^4,-K_{S_{2,2}^4}\otimes\mathcal I_p\bigr)
\to \left(\mathfrak m_p/\mathfrak m_p^2\right)\otimes(-K_{S_{2,2}^4})|_p\).
Both vector spaces have dimension \(2\). Its kernel is zero, since a
nonzero element in the kernel would define a member of
\(\mathcal C_p\) singular at \(p\). Therefore, this map is an
isomorphism. It follows that there is a member
\(C\in\mathcal C_p\) whose strict transform \(C'\) passes through
\(z\).

Suppose first that \(C\) is irreducible. By the convexity argument,
we may assume that \(C\not\subseteq\Supp(D)\). Then
\[
2-r
=
C'\cdot D'
\geq
\mult_z(D')
>
\frac2{\lambda}-r,
\]
which is impossible.

Therefore, we have \(C=L_1+L_2\), where \(p\in L_1\) and
\(p\notin L_2\). Put \(d\coloneqq L_1\cdot L_2\). Since
\((-K_{S_{2,2}^4})\cdot L_i=1\), we have
\(L_1^2=L_2^2=1-d\).

Let \(\nu\colon\widetilde S_{2,2}^4\to S_{2,2}^4\) be the
minimal resolution. The strict transform \(\overline L_i\) of \(L_i\)
is a smooth rational curve satisfying
\(-K_{\widetilde S_{2,2}^4}\cdot\overline L_i=1\). Adjunction gives
\(\overline L_i^2=-1\). If \(L_i\) does not pass through the
\(A_2\)-point, then \(L_i^2=-1\). If \(L_i\) passes through the
\(A_2\)-point, then, after relabelling the exceptional curves
\(Z_1\) and \(Z_2\), its numerical pullback is
\[
\nu^{\ast}L_i
=
\overline L_i+\frac23Z_1+\frac13Z_2,
\]
and hence \(L_i^2=-\frac13\). Thus,
\(d\in\left\{2,\frac43\right\}\).

Write \(D=b_1L_1+b_2L_2+\Omega_0\), where \(L_1,L_2\not\subseteq\Supp(\Omega_0)\), and put
\(\mu_0\coloneqq\min\{b_1,b_2\}\). We have \(\mu_0<1\). Indeed, if
\(\mu_0\geq1\), then \(D-C\) is effective and
\(\Q\)-linearly trivial, so \(D=C\), contradicting the fact that
\((S_{2,2}^4,\lambda C)\) is log canonical at \(p\).

Set
\[
D_{\mu_0}
\coloneqq
\frac{D-\mu_0C}{1-\mu_0}
\sim_{\Q}
-K_{S_{2,2}^4}.
\]
Write
\(\sigma^{\ast}D_{\mu_0}=D_{\mu_0}'+r_{\mu_0}F\). Since \(C\) is smooth
at \(p\), we have
\[
\begin{aligned}
\lambda D'+(\lambda r-1)F
=
(1-\mu_0)
\left(
\lambda D_{\mu_0}'
+(\lambda r_{\mu_0}-1)F
\right)+
\mu_0
\left(
\lambda C'+(\lambda-1)F
\right).
\end{aligned}
\]
The second pair on the right-hand side is log canonical at \(z\).
Consequently, the transformed pair associated with \(D_{\mu_0}\)
must remain non-log-canonical at \(z\); otherwise, convexity of log
canonicity would imply that the transformed pair associated with \(D\)
is log canonical at \(z\).

Replacing \(D\) by \(D_{\mu_0}\), and retaining the notation
\(D,D'\), and \(r\), we may assume that at least one of \(L_1,L_2\)
is not contained in \(\Supp(D)\), while
\[
\mult_z(D')+r>\frac2{\lambda}
\]
continues to hold. If \(L_1\not\subseteq\Supp(D)\), then, since
\(z\in L_1'\),
\[
1-r
=
L_1'\cdot D'
\geq
\mult_z(D')
>
\frac2{\lambda}-r,
\]
which is impossible. Therefore,
\(L_1\subseteq\Supp(D)\) and \(L_2\not\subseteq\Supp(D)\). We may thus
write \(D=aL_1+\Omega\), where \(a>0\) and \(L_1,L_2\not\subseteq\Supp(\Omega)\).

Since \(1=L_2\cdot D=ad+L_2\cdot\Omega\), we have \(ad\leq1\).
Moreover, \(L_1\cdot\Omega=1-aL_1^2=1+a(d-1)\).
In particular, \(a<1\), and hence \(\lambda a<1\).

Let \(L_1'\) and \(\Omega'\) be the strict transforms on
\(\widehat S_{2,2}^4\). Applying \cite[Lemma 2.5(iii)]{MG14} to \(L_1'\),
we obtain \(1<L_1'\cdot\left(\lambda\Omega'+(\lambda r-1)F
\right)=\lambda\left(L_1\cdot\Omega+a\right)-1=\lambda(1+ad)-1\leq2\lambda-1<\frac12\), which is impossible. Therefore, \(\alpha_p(S_{2,2}^4)\geq\frac34\).

For the surface \(S_{2,2}^5\), by \cite[Introduction, p. 1]{Den3}, one has \(\delta_q(S_{2,2}^5)\geq\frac{15}{7}\) at every smooth point \(q\in S_{2,2}^5\). Since \(\delta_q(S_{2,2}^5)\leq3\alpha_q(S_{2,2}^5)\), we obtain \(\alpha_q(S_{2,2}^5)\geq \frac13\delta_q(S_{2,2}^5) \geq \frac57\).
\end{proof}

\begin{notation}\label{not:standard-partial-resolution}
Let \(n\geq2\), \(m\geq2\), and \(k\geq0\). Let
\(\ell\in|\mathcal O_{\P(1,1,n)}(1)|\), and choose distinct smooth points \(a_1,\dots,a_m\in\ell\) and \(b_1,\dots,b_k\in\P(1,1,n)\setminus\ell\). We assume that the points \(b_1,\dots,b_k\) are general. In particular,
no two of them lie on the same member of
\(|\mathcal O_{\P(1,1,n)}(1)|\).

Let \(Y_{n,m}^k\) be the blow-up of \(\P(1,1,n)\) at
\(a_1,\dots,a_m,b_1,\dots,b_k\), and let \(L_{n,m}^k\) be the strict
transform of \(\ell\). We denote the contraction of \(L_{n,m}^k\) by \(\varphi_{n,m}^k\colon Y_{n,m}^k\to S_{n,m}^k\).
\end{notation}

\begin{lemma}\label{lem:elementary-transform-target}
Let \(n\geq3\), \(m\geq2\), and \(k\geq2\). With the notation of
Notation \ref{not:standard-partial-resolution}, for every
\(i\in\{1,\dots,k\}\), let
\(\ell_i\in|\mathcal O_{\P(1,1,n)}(1)|\) be the unique member containing
\(b_i\), and let \(R_i\subset Y_{n,m}^k\) be its strict transform. Then
\(\ell_i\neq\ell\), the curve \(R_i\) can be contracted, and its contraction \(\rho_i\colon Y_{n,m}^k\to Y^{-}\) satisfies \(Y^{-}\cong Y_{n-1,m}^{k-1}\). Moreover, if \(L^{-}\coloneqq(\rho_i)_{\ast}L_{n,m}^k\), then
\(L^{-}=L_{n-1,m}^{k-1}\), and the contraction of \(L^{-}\) gives
\(S_{n-1,m}^{k-1}\).

For every smooth point
\(y\in Y_{n,m}^k\setminus L_{n,m}^k\), one can choose \(i\) such that
\(y\notin R_i\). For such a choice, \(\rho_i(y)\notin L^{-}\).
\end{lemma}

\begin{proof}
Write \(\P(1,1,n)=\P(1,1,n)_{x,y,z}\). Since \(b_i\) is smooth, its first
two coordinates do not vanish simultaneously. Members of
\(|\mathcal O_{\P(1,1,n)}(1)|\) have equations
\(\alpha x+\beta y=0\). Hence, there is a unique member \(\ell_i\)
containing \(b_i\). Since \(b_i\notin\ell\), we have
\(\ell_i\neq\ell\).

Let \(H\) be the pullback of \(\mathcal O_{\P(1,1,n)}(1)\), and let
\(F_i\) be the exceptional curve over \(b_i\). By the generality of the
points \(b_1,\dots,b_k\), no other blow-up center lies on \(\ell_i\).
Consequently,
\[
    R_i\equiv H-F_i,
    \qquad
    R_i^2=\frac1n-1=-\frac{n-1}{n}.
\]

We describe the contraction using the minimal resolution of
\(\P(1,1,n)\). Let
\(\mu_n\colon\mathbb F_n\to\P(1,1,n)\) be the contraction of the negative
section \(C_n\subset\mathbb F_n\), where \(C_n^2=-n\). The curve
\(\ell_i\) is the image of a fibre \(f_i\) of \(\mathbb F_n\). The point
\(b_i\) lies on \(f_i\setminus C_n\). After blowing up \(b_i\), the strict
transform of \(f_i\) is a \((-1)\)-curve. Contracting it is the standard
elementary transformation \(\mathbb F_n\dashrightarrow\mathbb F_{n-1}\).
Indeed, the image of \(C_n\) has self-intersection
\(-n+1=-(n-1)\). Contracting this negative section gives
\(\P(1,1,n-1)\).

The points \(a_1,\dots,a_m\) lie on the member \(\ell\), which is distinct
from \(\ell_i\). The points \(b_j\) with \(j\neq i\) do not lie on
\(\ell_i\). Hence, all remaining blow-ups commute with this elementary
transformation. It follows that \(Y^{-}\cong Y_{n-1,m}^{k-1}\).

The image of \(\ell\) under the elementary transformation is a member of
\(|\mathcal O_{\P(1,1,n-1)}(1)|\). Equivalently, one can verify its
self-intersection directly. Since
\[
    \left(L_{n,m}^k\right)^2=\frac1n-m,
    \qquad
    L_{n,m}^k\cdot R_i=\frac1n,
\]
we obtain
\[
    \left((\rho_i)_{\ast}L_{n,m}^k\right)^2
    =
    \frac1n-m
    -
    \frac{\left(\frac1n\right)^2}{-\frac{n-1}{n}}
    =
    \frac1{n-1}-m.
\]
The images of \(a_1,\dots,a_m\) remain on this member, while the images of
\(b_j\), \(j\neq i\), remain general smooth  points outside it. Hence, we have \((\rho_i)_{\ast}L_{n,m}^k=L_{n-1,m}^{k-1}\).

Distinct members of
\(|\mathcal O_{\P(1,1,n)}(1)|\) meet only at the quotient point
\([0:0:1]\). Hence, the curves \(R_1,\dots,R_k\) meet one another only at
that quotient point. A smooth point of \(Y_{n,m}^k\) therefore belongs to
at most one of the curves \(R_i\). Since \(k\geq2\), for every smooth
\(y\in Y_{n,m}^k\setminus L_{n,m}^k\), one can choose \(i\) such that
\(y\notin R_i\).

Set-theoretically, we have \(\rho_i^{-1}(L^{-})=L_{n,m}^k\cup R_i\).
Hence, \(y\notin L_{n,m}^k\cup R_i\) implies
\(\rho_i(y)\notin L^{-}\).
\end{proof}

\begin{lemma}\label{lem:smooth-alpha-reduction}
Let \(n\geq3\), \(m\geq2\), and \(k\geq2\). Consider the contractions:
\[
    \varphi\colon Y_{n,m}^k\longrightarrow S_{n,m}^k
    \qquad\text{and}\qquad
    \varphi^{-}\colon Y_{n-1,m}^{k-1}\longrightarrow
    S_{n-1,m}^{k-1}.
\]
Let \(\eta\in\Q_{>0}\). Suppose that \(\alpha_q(S_{n-1,m}^{k-1})\geq\eta\)
for every smooth point \(q\in S_{n-1,m}^{k-1}\). Then we have \(\alpha_p(S_{n,m}^k)\geq\eta\) for every smooth point \(p\in S_{n,m}^k\).
\end{lemma}

\begin{proof}
Put \(Y\coloneqq Y_{n,m}^k\), \(X\coloneqq S_{n,m}^k\), and
\(L\coloneqq L_{n,m}^k\). We first compute the discrepancy along \(L\).
Since
\[
    L^2=\frac{1-mn}{n}
    \qquad\text{and}\qquad
    K_Y\cdot L=m-\frac{n+2}{n},
\]
we have \(K_Y = \varphi^{\ast}K_X-\gamma_{n,m}L\), where \(\gamma_{n,m} \coloneqq \frac{mn-n-2}{mn-1}\). Note that \(\gamma_{n,m}\geq0\).

Fix a rational number \(\lambda\) with \(0<\lambda<\eta\). Let
\(D_X\sim_{\Q}-K_X\) be an effective \(\Q\)-divisor. Suppose that
\((X,\lambda D_X)\) is not log canonical at a smooth point \(p\in X\).
Write \(\varphi^{\ast}D_X=D_Y+\nu L\), where \(D_Y\geq0\) and \(L\not\subseteq\Supp(D_Y)\). Put \(G \coloneqq D_Y+\left(\nu+\gamma_{n,m}\right)L\).
Then \(G\geq0\) and \(G\sim_{\Q}-K_Y\).

Since \(p\) is smooth, \(p\neq\varphi(L)\). Hence,
\(\varphi^{-1}(p)\) consists of a single smooth point
\(y\in Y\setminus L\), and \(\varphi\) is an isomorphism near \(y\).
The divisors \(G\) and \(D_Y\) coincide near \(y\). Consequently,
\((Y,\lambda G)\) is not log canonical at \(y\).

By Lemma \ref{lem:elementary-transform-target}, choose \(i\) such that
\(y\notin R_i\). Let \(\rho_i\colon Y\to Y^{-}\cong Y_{n-1,m}^{k-1}\) be the contraction of \(R_i\), and put \(L^{-}\coloneqq(\rho_i)_{\ast}L\) and \(G^{-}\coloneqq(\rho_i)_{\ast}G\).
The morphism \(\rho_i\) is an isomorphism near \(y\). Moreover,
\(\rho_i(y)\notin L^{-}\). Choosing compatible canonical divisors, we
have \((\rho_i)_{\ast}K_Y=K_{Y^{-}}\), and therefore \(G^{-}\sim_{\Q}-K_{Y^{-}}\).
The pair \((Y^{-},\lambda G^{-})\) is not log canonical at the smooth
point \(\rho_i(y)\).

Let \(X^{-}\coloneqq S_{n-1,m}^{k-1}\) and \(B^{-}\coloneqq(\varphi^{-})_{\ast}G^{-}\).
Since \(\rho_i(y)\notin L^{-}\), the morphism \(\varphi^{-}\) is an
isomorphism near \(\rho_i(y)\). Moreover, we have \(0\leq B^{-}\sim_{\Q}-K_{X^{-}}\).
Hence, the pair \((X^{-},\lambda B^{-})\) is not log canonical at the
smooth point
\(\varphi^{-}(\rho_i(y))\). This contradicts the assumption that
\(\alpha_q(X^{-})\geq\eta>\lambda\) at every smooth point \(q\in X^{-}\).

Since \(\lambda<\eta\) was arbitrary, we conclude that
\(\alpha_p(X)\geq\eta\) for every smooth point \(p\in X\).
\end{proof}

\subsubsection{Case: \(S_{3,2}^6\)}

Let \(\pi\colon S\to \P(1,1,3)\) be the blow-up of \(\P(1,1,3)\) at seven general smooth  points \(p_1,\dots,p_7\). Note that the surface \(S\) is isomorphic to a complete intersection of two hypersurfaces of degree \(4\) in \(\P(1,1,2,2,3)\) by Proposition \ref{prop:ci hyp}. Let \(L\) be the strict transform of the curve \(\ell\in |\mathcal{O}_{\P(1,1,3)}(1)|\) passing through the point \(p_1\). Let \(\pi_1\colon S_1\to S\) be the blow-up of a general smooth point \(r\in L\) such that \(\pi(r)\neq p_1\). Then, by contracting the strict transform \(L_1\) of \(L\), we obtain the birational morphism \(\varphi\colon S_1\to S_{3,2}^6\). For the singular point \(o\in S_1\) of type \(\frac{1}{3}(1,1)\), let \(\pi_2\colon S_2\to S_1\) be the weighted blow-up with weights \((1,1)\). We consider the following diagram.
\begin{figure}[H]
  \centering
  \begin{tikzpicture}
  [-,auto,node distance=1.5cm, thick,main node/.style={circle,draw,font=\sffamily \Large\bfseries}]
    \node[text=black] (1) {$\mathbb{P}(1,1,3)$};
    \node[text=black] (2) [right=1.5cm of 1] {$S$};
    \node[text=black] (3) [right of=2, above of=2] {$S_1$};
    \node[text=black] (4) [above of=3] {$S_2$};
    \node[text=black] (5) [right of=3, below of=3] {$S_{3,2}^6$};
   
    \path[every node/.style={font=\sffamily}]
      (2) edge[->] node [above, pos=0.4] {\(\pi\)} (1)
      (3) edge[->] node [above, left=0.5mm, pos=0.3] {\(\pi_1\)} (2)
      (4) edge[->] node [left] {\(\pi_2\)} (3)
      (3) edge[->] node [above, right=0.5mm, pos=0.3] {\(\varphi\)} (5);
  \end{tikzpicture}
    \caption{}\label{fig:326}
\end{figure}

\begin{lemma}\label{lem;326 sing}
    Let \(p\in S_{3,2}^6\) be the singular point of type \(\frac{1}{5}(1,3)\). Then we have \(\delta_p(S_{3,2}^6)>1\).
\end{lemma}
\begin{proof}
    The linear system \(|\mathcal O_{\P(1,1,3)}(4)|\) has projective dimension \(6\). Hence, there is a unique member \(c\) passing through \(p_2,\dots,p_7\). For a general configuration, this member is irreducible and has the transverse intersection properties used below. Let \(C\) be its strict transform under \(\pi\). By the generality assumption, \(r\notin C\). Then we have the following:
\begin{align*}
    -K_S&\equiv L+C,\quad L^2=C^2=-\frac{2}{3},\quad L\cdot C=\frac{4}{3},\\
    -K_{S_1}&\equiv L_1+C_1,\quad L_1^2=-\frac{5}{3},\quad C_1^2=-\frac{2}{3}, \quad L_1\cdot C_1=\frac{4}{3},
\end{align*}
where \(L_1\) and \(C_1\) are the strict transforms of \(L\) and \(C\), respectively. We also have \(-K_{S_{3,2}^6}\equiv R_S\) and \((-K_{S_{3,2}^6})^2=\frac{2}{5}\) such that \(\varphi_{\ast}(C_1)=R_S\). We obtain \(\varphi^{\ast}(R_S)=C_1+\frac{4}{5}L_1\).
Since \(C_1\) is a negative curve, the pseudoeffective threshold of \(L_1\) is \(4/5\). 
The Zariski decomposition $\varphi^{\ast}(-K_{S_{3,2}^6})-tL_1=P(t)+N(t)$ is given by
\[
P(t) = 
\begin{dcases}
    \varphi^{\ast}(-K_{S^6_{3,2}})-tL_1,& t\in \left[0,\frac{3}{10}\right],\\
    \left(\frac{4}{5}-t\right)(2C_1+L_1),& t\in \left[\frac{3}{10},\frac{4}{5}\right],
\end{dcases}
\quad
N(t) = 
\begin{dcases}
    0,& t\in \left[0,\frac{3}{10}\right],\\
    \left(2t-\frac{3}{5}\right)C_1,& t\in \left[\frac{3}{10},\frac{4}{5}\right].
\end{dcases}
\]
Then \(P(t)^2\) is given by
\[
    P(t)^2=
    \begin{dcases}
        \frac{2}{5}-\dfrac{5}{3}t^2, & t\in \left[0,\frac{3}{10}\right],\\
        \left(\frac{4}{5}-t\right)^2, & t\in \left[\frac{3}{10},\frac{4}{5}\right].
    \end{dcases}
\]
Hence, the \(S\)-invariant of \(L_1\) is \(S(L_1)=\frac{5}{2}\left(\frac{21}{200}+\frac{1}{24}\right)=\frac{11}{30}\). Since \(K_{S_1}=\varphi^{\ast}(K_{S_{3,2}^6})-\frac{1}{5}L_1\), we have \(A_{S_{3,2}^6}(L_1)=\frac{4}{5}\). 
Therefore,
\[
\frac{A_{S_{3,2}^6}(L_1)}{S_{S_{3,2}^6}(L_1)} = \frac{24}{11}>1.
\] 
Furthermore, \(P(t)\cdot L_1\) is given by
\[
    P(t)\cdot L_1=
    \begin{dcases}
        \frac{5}{3}t, & t\in \left[0,\frac{3}{10}\right],\\
        \frac{4}{5}-t, & t\in \left[\frac{3}{10},\frac{4}{5}\right].
    \end{dcases}
\]
For a smooth point \(q\in L_1\setminus C_1\), the function \(h(t)\) in
\eqref{function h(t)} is given by
\[
    h(t)=
    \begin{dcases}
        \frac{25}{18}t^2, & t\in \left[0,\frac{3}{10}\right],\\
        \frac{1}{2}\left(\dfrac{4}{5}-t\right)^2, &  t\in \left[\frac{3}{10},\frac{4}{5}\right].
    \end{dcases}
\]
Thus, \(S(W_{\bullet,\bullet}^{L_1};q)=\frac{1}{6}\). Since \(A_{L_1,\Phi_{L_1}}(q)=1\), we obtain
\[
\frac{A_{L_1,\Phi_{L_1}}(q)}
{S(W_{\bullet,\bullet}^{L_1};q)}
=6>1.
\]
The curves \(C_1\) and \(L_1\) meet at the singular point of type
\(\frac{1}{3}(1,1)\) and at one additional smooth point.
Denote this smooth intersection point by \(q_0\). For $t\in [\frac{3}{10},\frac{4}{5}]$ the negative part $N(t)$ contains the component \(\left(2t-\frac{3}{5}\right)C_1\).
Hence, for $t\in [\frac{3}{10},\frac{4}{5}]$, we have
\[
(P(t)\cdot L_1)(N(t)\cdot L_1)_{q_0}
=
\left(\frac{4}{5}-t\right)\left(2t-\frac{3}{5}\right).
\]
For the smooth point \(q_0\), the function \(h(t)\) in \eqref{function h(t)} is given by
\[
    h(t)=
    \begin{dcases}
        \frac{25}{18}t^2, &  t\in \left[0,\frac{3}{10}\right],\\
        \left(\frac{4}{5}-t\right)\left(\dfrac{3}{2}t-\dfrac{1}{5}\right),
        &  t\in \left[\frac{3}{10},\frac{4}{5}\right].
    \end{dcases}
\]
It follows that \(S(W_{\bullet,\bullet}^{L_1};q_0)=\frac{3}{8}\).
Since \(A_{L_1,\Phi_{L_1}}(q_0)=1\), we obtain
\[
\frac{A_{L_1,\Phi_{L_1}}(q_0)}
{S(W_{\bullet,\bullet}^{L_1};q_0)}
=
\frac{8}{3}>1.
\]

The preceding calculations will be used below to control the points on
the strict transform of \(L_1\) away from the exceptional divisor of
\(\pi_2\). We next compute the corresponding data along the divisor
\(E\). Set \(\psi\coloneqq\varphi\circ\pi_2\). Then we obtain 
\begin{align*}
    \psi^{\ast}(R_S)=\pi_2^{\ast}\left(C_1+\frac{4}{5}L_1\right)=C_2+\frac{4}{5}L_2+\frac{3}{5}E,
\end{align*}
where \(E\) is the \(\pi_2\)-exceptional divisor, \(C_2\) and \(L_2\) are the strict transforms of \(C_1\) and \(L_1\), respectively. The intersection numbers are the following:
\begin{align*}
    L_2^2=-2,\quad C_2^2=-1,\quad L_2\cdot C_2=1,\quad E^2=-3, \quad L_2\cdot E=C_2\cdot E=1.
\end{align*}
Since the intersection matrix corresponding to the curves \(C_2\) and \(L_2\) is negative definite, the pseudoeffective threshold of \(E\) is \(\frac{3}{5}\).
The Zariski decomposition $\psi^{\ast}(-K_{S_{3,2}^6})-tE=P(t)+N(t)$
is given by
\[
    P(t)=
    \begin{dcases}
        C_2+\left(\frac{4}{5}-\frac{1}{2}t\right)L_2+\left(\frac{3}{5}-t\right)E, & t\in\left[0,\frac{4}{15}\right],\\
        \left(\frac{3}{5}-t\right)(3C_2+2L_2+E), & t\in\left[\frac{4}{15},\frac35\right],
    \end{dcases}
\]
and
\[
N(t)=
    \begin{dcases}
        \frac{1}{2}tL_2, & t\in\left[0,\frac{4}{15}\right],\\
        \left(-\frac{4}{5}+3t\right)C_2+\left(-\frac{2}{5}+2t\right)L_2, & t\in\left[\frac{4}{15},\frac35\right].
    \end{dcases}
\]
Then \(P(t)^2\) is given by
\[
    P(t)^2=
    \begin{dcases}
        -\dfrac{1}{10}(5t+2)(5t-2), & t\in \left[0,\frac{4}{15}\right],\\
        2\left(\dfrac{3}{5}-t\right)^2, & t\in \left[\frac{4}{15},\frac{3}{5}\right].
    \end{dcases}
\]
Hence, the \(S\)-invariant of \(E\) is \(S(E)=\frac{5}{2}\left(\frac{184}{2025}+\frac{2}{81}\right)=\frac{13}{45}\). Since \(K_{S_2}=\psi^{\ast}(K_{S_{3,2}^6})-\frac{1}{5}L_2-\frac{2}{5}E\), we have \(A_{S_{3,2}^6}(E)=\frac{3}{5}\). Therefore, 
\[
    \frac{A_{S_{3,2}^6}(E)}{S_{S_{3,2}^6}(E)} = \frac{27}{13}>1.
\] 
Furthermore, \(P(t)\cdot E\) is given by
\[
    P(t)\cdot E=
    \begin{dcases}
        \dfrac{5}{2}t, & t\in \left[0,\frac{4}{15}\right],\\
        \dfrac{6}{5}-2t, & t\in \left[\frac{4}{15},\frac{3}{5}\right].
    \end{dcases}
\]

For the point \(q_1\in C_2\cap E\), the function \(h(t)\) in \eqref{function h(t)} is given by
\begin{align*}
    h(t)=
    \begin{dcases}
        \frac{25}{8}t^2, & t\in \left[0,\frac{4}{15}\right],\\
        \left(\frac{6}{5}-2t\right)\left(-\frac{4}{5}+3t\right)+\frac{1}{2}\left(\frac{6}{5}-2t\right)^2, & t\in \left[\frac{4}{15},\frac{3}{5}\right].
    \end{dcases}
\end{align*}
Hence, \(S(W_{\bullet,\bullet}^E;q_1) = 5\left(\frac{8}{405}+\frac{5}{81}\right) = \frac{11}{27}\).
Since no other exceptional component contributes to the different at the point \(C_2 \cap E\), we have \(A_{E,\Phi_E}(q_1)=1\).
Thus,  
\[
\frac{A_{E,\Phi_E}(q_1)}
{S(W_{\bullet,\bullet}^E;q_1)}
=
\frac{27}{11}>1.
\]
For the point \(q_2\in L_2\cap E\), the function \(h(t)\) in \eqref{function h(t)} is given by 
\begin{align*}
    h(t)=
    \begin{dcases}
        \frac{5}{4}t^2+\frac{25}{8}t^2, & t\in \left[0,\frac{4}{15}\right],\\
        \left(\frac{6}{5}-2t\right)\left(-\frac{2}{5}+2t\right)+\frac{1}{2}\left(\frac{6}{5}-2t\right)^2, & t\in \left[\frac{4}{15},\frac{3}{5}\right].
    \end{dcases}
\end{align*}
Thus, \(S(W_{\bullet,\bullet}^E;q_2) = 5\left(\frac{56}{2025}+\frac{26}{405}\right) = \frac{62}{135}\).
By Remark \ref{rem:AZ-higher-model}, the divisor \(L_2\) contributes to the
different on \(E\). Hence, \(A_{E,\Phi_E}(q_2)=A_{S_{3,2}^6}(L_2)=\frac{4}{5}\). Then we have
\[
\frac{A_{E,\Phi_E}(q_2)}
{S(W_{\bullet,\bullet}^E;q_2)}
=
\frac{54}{31}>1.
\]
For \(q_3\in E\setminus(C_2\cup L_2)\), we have \(A_{E,\Phi_E}(q_3)=1\) and
\begin{align*}
    h(t)=
    \begin{dcases}
        \frac{25}{8}t^2, & t\in \left[0,\frac{4}{15}\right],\\
        \frac{1}{2}\left(\frac{6}{5}-2t\right)^2, & t\in \left[\frac{4}{15},\frac{3}{5}\right].
    \end{dcases}
\end{align*}
Then \(S(W_{\bullet,\bullet}^E;q_3) = 5\left(\frac{8}{405}+\frac{2}{81}\right) = \frac{2}{9}\).
Thus, we obtain
\[
\frac{A_{E,\Phi_E}(q_3)}
{S(W_{\bullet,\bullet}^E;q_3)}
=
\frac{9}{2}>1.
\]

The canonical divisor formula gives \(K_{S_2}+B_{S_2}=\psi^{\ast}K_{S_{3,2}^6}\), where \(B_{S_2}=\frac15L_2+\frac25E\).
Applying Lemma \ref{lem:AZ-crepant-model} with the divisor
\(E\), the preceding calculations give
\[
\inf_{q\in E}
\delta_q(S_2,B_{S_2};\psi^{\ast}(-K_{S_{3,2}^6}))
\geq
\min\left\{
    \frac{27}{13},
    \frac{27}{11},
    \frac{54}{31},
    \frac92
\right\}
=
\frac{54}{31}.
\]

The divisors \(L_1\subset S_1\) and \(L_2\subset S_2\) define the same
divisorial valuation over \(S_{3,2}^6\). Since \(\pi_2\) is an
isomorphism over \(L_2\setminus E\), pullback identifies the restricted
graded linear series and their induced filtrations at corresponding
points of \(L_1\) and \(L_2\setminus E\). Hence, the calculations for
\(L_1\) give
\[
\inf_{q\in L_2\setminus E}
\delta_q(S_2,B_{S_2};\psi^{\ast}(-K_{S_{3,2}^6}))
\geq
\min\left\{
    \frac{24}{11},
    6,
    \frac83
\right\}
=
\frac{24}{11}.
\]

Since \(\psi^{-1}(p)=L_2\cup E\), Lemma \ref{lem:AZ-crepant-model} yields
\[
    \delta_p(S_{3,2}^6)
    \geq
    \min\left\{
        \frac{24}{11},
        \frac{54}{31}
    \right\}
    =
    \frac{54}{31}
    >
    1.\qedhere
\]
\end{proof}

Now, we show that \(\delta_p(S_{3,2}^6)>1\) for smooth points \(p\).

\begin{lemma}\label{lem;326 smooth}
Let \(p\in S_{3,2}^6\) be a smooth point. Then we have \(\alpha_p(S_{3,2}^6)\geq\frac{7}{10}\). In particular, \(\delta_p(S_{3,2}^6)>1\).
\end{lemma}

\begin{proof}
Apply Lemma \ref{lem:smooth-alpha-reduction} with
\((n,m,k)=(3,2,6)\). The target of the elementary transformation is
\(S_{2,2}^5\). By Lemma \ref{lem:base-22-smooth-alpha}, we obtain \(\alpha_q(S_{2,2}^5)\geq\frac57>\frac7{10}\) for every smooth point \(q\in S_{2,2}^5\). Hence, \(\alpha_p(S_{3,2}^6)\geq\frac{7}{10}\).

By Theorem \ref{thm:alpha delta}, we have
\[
    \delta_p(S_{3,2}^6)
    \geq
    \frac32\alpha_p(S_{3,2}^6)
    \geq
    \frac{21}{20}
    >
    1.\qedhere
\]
\end{proof}

\begin{theorem}\label{thm:326 stable}
    The surface \(S_{3,2}^6\) is K-stable.
\end{theorem}
\begin{proof}
Lemmas \ref{lem;326 sing} and \ref{lem;326 smooth} give \(\delta(S_{3,2}^6)\geq\frac{21}{20}>1\).
Therefore, Theorems \ref{thm:delta} and
\ref{thm:K-stable-uniform} imply that \(S_{3,2}^6\) is K-stable.
\end{proof}

\subsubsection{Case: \(S_{4,2}^7\)}

Let \(\pi\colon S'\to \P(1,1,4)\) be the blow-up of \(\P(1,1,4)\) at eight general smooth  points \(p_1,\dots,p_8\). Note that the surface \(S'\) is isomorphic to a degree \(6\) hypersurface embedded in \(\P(1,1,2,3)\) by Proposition \ref{prop:ci hyp}. Let \(L\) be the strict transform of the curve \(\ell\in |\mathcal{O}_{\P(1,1,4)}(1)|\) passing through the point \(p_1\). Let \(\pi_1\colon S'_1\to S'\) be the blow-up of a general smooth point \(r\in L\) such that \(\pi(r)\neq p_1\). Then, by contracting the strict transform \(L_1\) of \(L\), we obtain the birational morphism \(\varphi\colon S'_1\to S_{4,2}^7\). For the singular point \(o\in S'_1\) of type \(\frac{1}{4}(1,1)\), let
\(\pi_2\colon S'_2\to S'_1\) be the weighted blow-up with weights \((1,1)\). We consider the following diagram.
\begin{figure}[H]
  \centering
  \begin{tikzpicture}
  [-,auto,node distance=1.5cm, thick,main node/.style={circle,draw,font=\sffamily \Large\bfseries}]
    \node[text=black] (1) {$\mathbb{P}(1,1,4)$};
    \node[text=black] (2) [right=1.5cm of 1] {$S'$};
    \node[text=black] (3) [right of=2, above of=2] {$S'_1$};
    \node[text=black] (4) [above of=3] {$S'_2$};
    \node[text=black] (5) [right of=3, below of=3] {$S_{4,2}^7$};
   
    \path[every node/.style={font=\sffamily}]
      (2) edge[->] node [above, pos=0.4] {\(\pi\)} (1)
      (3) edge[->] node [above, left=0.5mm, pos=0.3] {\(\pi_1\)} (2)
      (4) edge[->] node [left] {\(\pi_2\)} (3)
      (3) edge[->] node [above, right=0.5mm, pos=0.3] {\(\varphi\)} (5);
  \end{tikzpicture}
    \caption{}\label{fig:427}
\end{figure}

\begin{lemma}\label{lem;427 sing}
    Let \(p\in S_{4,2}^7\) be the singular point of type \(\frac{1}{7}(1,4)\). Then we have \(\delta_p(S_{4,2}^7)>1\).
\end{lemma}
\begin{proof}
    The linear system \(|\mathcal O_{\P(1,1,4)}(5)|\) has projective dimension \(7\). Hence, there is a unique member \(c\) passing through \(p_2,\dots,p_8\). For a general configuration, this member is irreducible and has the transverse intersection properties used below. Let \(C\) be its strict transform under \(\pi\). By the generality assumption, \(r\notin C\). Then we have the following:
\begin{align*}
    -K_{S'}&\equiv L+C,\quad L^2=C^2=-\frac{3}{4},\quad L\cdot C=\frac{5}{4},\\
    -K_{S'_1}&\equiv L_1+C_1,\quad L_1^2=-\frac{7}{4},\quad C_1^2=-\frac{3}{4}, \quad L_1\cdot C_1=\frac{5}{4},
\end{align*}
where \(L_1\) and \(C_1\) are the strict transforms of \(L\) and \(C\), respectively. Then we have \(-K_{S_{4,2}^7}\equiv R_{S'}\) and \( (-K_{S_{4,2}^7})^2=\frac{1}{7}\) such that \(\varphi_{\ast}(C_1)=R_{S'}\). We obtain \(\varphi^{\ast}(R_{S'})=C_1+\frac{5}{7}L_1\).

Since \(C_1\) is a negative curve, the pseudoeffective threshold of \(L_1\) is \(\frac{5}{7}\). The positive and negative parts of the Zariski decomposition of \(\varphi^{\ast}(-K_{S_{4,2}^7})-tL_1\) are
\[
P(t)=
\begin{dcases}
    \varphi^{\ast}(-K_{S_{4,2}^7})-tL_1,&  t\in\left[0,\frac{4}{35}\right],\\
    \left(\frac{5}{7}-t\right)\left(\frac{5}{3}C_1+L_1\right),&  t\in\left[\frac{4}{35},\frac{5}{7}\right],
\end{dcases}
\quad
N(t)=
\begin{dcases}
    0, & t\in\left[0,\frac{4}{35}\right],\\
    \frac{1}{3}\left(-\frac{4}{7}+5t\right)C_1, & t\in\left[\frac{4}{35},\frac{5}{7}\right].
\end{dcases}
\]
Then we obtain
\begin{align*}
    P(t)^2=
    \begin{dcases}
        \frac{1}{7}-\frac{7}{4}t^2, & t\in\left[0,\frac{4}{35}\right],\\
        \frac{1}{3}\left(\frac{5}{7}-t\right)^2, & t\in\left[\frac{4}{35},\frac{5}{7}\right].
    \end{dcases}
\end{align*}
Hence, the \(S\)-invariant of \(L_1\) is \(S(L_1)=7\left(\frac{284}{18375}+\frac{3}{125}\right)=\frac{29}{105}\). Also, since \(K_{S'_1}=\varphi^{\ast}(K_{S_{4,2}^7})-\frac{2}{7}L_1\), \(A_{S_{4,2}^7}(L_1)=\frac{5}{7}\). Therefore, 
\[
\frac{A_{S_{4,2}^7}(L_1)} {S_{S_{4,2}^7}(L_1)} = \frac{75}{29}>1.
\]
Moreover, we have
\begin{align*}
    P(t)\cdot L_1=
    \begin{dcases}
        \frac{7}{4}t, & t\in\left[0,\frac{4}{35}\right],\\
        \frac{1}{3}\left(\frac{5}{7}-t\right), & t\in\left[\frac{4}{35},\frac{5}{7}\right].
    \end{dcases}
\end{align*}
For a smooth point \(q\in L_1\setminus C_1\), we have
\[
h(t)=\frac{1}{2}(P(t)\cdot L_1)^2=
\begin{dcases}
\frac{49}{32}t^2,
& t\in\left[0,\frac{4}{35}\right],\\
\frac{1}{18}\left(\frac{5}{7}-t\right)^2,
& t\in\left[\frac{4}{35},\frac{5}{7}\right].
\end{dcases}
\]
Therefore, \(S(W_{\bullet,\bullet}^{L_1};q)=\frac{1}{15}\). Since \(A_{L_1,\Phi_{L_1}}(q)=1\), we obtain 
\[
\frac{A_{L_1,\Phi_{L_1}}(q)}
{S(W_{\bullet,\bullet}^{L_1};q)}
=
15>1.
\]

The curves \(C_1\) and \(L_1\) meet at the singular point of type \(\frac{1}{4}(1,1)\) and at one additional smooth point. Now, let \(q\in C_1\cap L_1\) be the smooth intersection point. For \(\frac{4}{35}\leq t\leq \frac{5}{7}\), the negative part is \(\frac{1}{3}\left(-\frac{4}{7}+5t\right)C_1\).
Hence, the local intersection term contributes
\[
(P(t)\cdot L_1)(N(t)\cdot L_1)_q
=
\frac{1}{3}\left(\frac{5}{7}-t\right)
\cdot
\frac{1}{3}\left(-\frac{4}{7}+5t\right).
\]
Therefore,
\[
h(t)=
\begin{dcases}
\frac{49}{32}t^2,
& t\in\left[0,\frac{4}{35}\right],\\
\frac{1}{18}\left(\frac{5}{7}-t\right)^2+
\frac{1}{3}\left(\frac{5}{7}-t\right)
\cdot
\frac{1}{3}\left(-\frac{4}{7}+5t\right),
& t\in\left[\frac{4}{35},\frac{5}{7}\right].
\end{dcases}
\]
It follows that \(S(W_{\bullet,\bullet}^{L_1};q)=\frac{26}{75}\).
Since \(A_{L_1,\Phi_{L_1}}(q)=1\), we obtain
\[
\frac{A_{L_1,\Phi_{L_1}}(q)}
{S(W_{\bullet,\bullet}^{L_1};q)}
=
\frac{75}{26}>1.
\]

The preceding calculations will be used below to control the points on
the strict transform of \(L_1\) away from the exceptional divisor of
\(\pi_2\). We next compute the corresponding data along the divisor
\(E\). Set \(\psi\coloneqq\varphi\circ\pi_2\). Then we obtain 
\begin{align*}
    \psi^{\ast}(R_{S'}) = \pi_2^{\ast}\left(C_1+\frac{5}{7}L_1\right) = C_2+\frac{5}{7}L_2+\frac{3}{7}E,
\end{align*}
where \(E\) is the \(\pi_2\)-exceptional divisor, \(C_2\) and \(L_2\) are the strict transforms of \(C_1\) and \(L_1\), respectively. The intersection numbers are the following:
\begin{align*}
    L_2^2=-2,\quad C_2^2=-1,\quad L_2\cdot C_2=1,\quad E^2=-4, \quad L_2\cdot E=C_2\cdot E=1.
\end{align*}
Since the intersection matrix corresponding to the curves \(C_2\) and \(L_2\) is negative definite, the pseudoeffective threshold of \(E\) is \(\frac{3}{7}\).

The positive and negative parts of the Zariski decomposition of \(\psi^{\ast}(-K_{S_{4,2}^7})-tE\) are
\begin{align*}
    P(t)=
    \begin{dcases}
        C_2+\left(\frac{5}{7}-\frac{1}{2}t\right)L_2+\left(\frac{3}{7}-t\right)E, & t\in \left[0,\frac{2}{21}\right],\\
        \left(\frac{3}{7}-t\right)(3C_2+2L_2+E), & t\in\left[\frac{2}{21},\frac{3}{7}\right],
    \end{dcases}
\end{align*}
and
\begin{align*}
    N(t)=
    \begin{dcases}
        \frac{1}{2}tL_2, & t\in \left[0,\frac{2}{21}\right],\\
        \left(-\frac{2}{7}+3t\right)C_2+\left(-\frac{1}{7}+2t\right)L_2, & t\in\left[\frac{2}{21},\frac{3}{7}\right].
    \end{dcases}
\end{align*}
We obtain 
\begin{align*}
    P(t)^2=
    \begin{dcases}
        \frac{1}{14}(2-49t^2), & t\in \left[0,\frac{2}{21}\right],\\
        \left(\frac{3}{7}-t\right)^2, & t\in\left[\frac{2}{21},\frac{3}{7}\right].
    \end{dcases}
\end{align*}
Hence, the \(S\)-invariant of \(E\) is \(S(E)=7\left(\frac{50}{3969}+\frac{1}{81}\right)=\frac{11}{63}\). Also, we have \(K_{S'_2}=\psi^{\ast}(K_{S_{4,2}^7})-\frac{2}{7}L_2-\frac{4}{7}E\), which implies that \(A_{S_{4,2}^7}(E)=\frac{3}{7}\). Therefore, 
\[
\frac{A_{S_{4,2}^7}(E)}{S_{S_{4,2}^7}(E)} = \frac{27}{11}>1.
\]
Moreover, we have
\begin{align*}
    P(t)\cdot E=
    \begin{dcases}
        \frac{7}{2}t, & t\in \left[0,\frac{2}{21}\right],\\
        \frac{3}{7}-t, & t\in\left[\frac{2}{21},\frac{3}{7}\right].
    \end{dcases}
\end{align*}

For the point \(q\in C_2\cap E\), we have 
\begin{align*}
    h(t)=
    \begin{dcases}
        \frac{49}{8}t^2, & t\in \left[0,\frac{2}{21}\right],\\
        \left(\frac{3}{7}-t\right)\left(-\frac{2}{7}+3t\right)+\frac{1}{2}\left(\frac{3}{7}-t\right)^2, & t\in\left[\frac{2}{21},\frac{3}{7}\right].
    \end{dcases}
\end{align*}
Hence, we have \(S(W_{\bullet,\bullet}^E;q)=\frac{10}{27}\). Since no other exceptional component contributes to the different at the point \(E\cap C_2\), we have \(A_{E,\Phi_E}(q)=1\).
Thus, we obtain
\[
\frac{A_{E,\Phi_E}(q)}
{S(W_{\bullet,\bullet}^E;q)}
=
\frac{27}{10}>1.
\]

For the point \(q\in L_2\cap E\), we have 
\begin{align*}
    h(t)=
    \begin{dcases}
        \frac{7}{4}t^2+\frac{49}{8}t^2, & t\in \left[0,\frac{2}{21}\right],\\
        \left(\frac{3}{7}-t\right)\left(-\frac{1}{7}+2t\right)+\frac{1}{2}\left(\frac{3}{7}-t\right)^2, & t\in\left[\frac{2}{21},\frac{3}{7}\right].
    \end{dcases}
\end{align*}
Then \(S(W_{\bullet,\bullet}^E;q) = 14\left(\frac{1}{441}+\frac{4}{189}\right) = \frac{62}{189}\). By Remark \ref{rem:AZ-higher-model}, the divisor \(L_2\) contributes to the different on \(E\). Hence, \(A_{E,\Phi_E}(q)=A_{S_{4,2}^7}(L_2)=\frac{5}{7}\). Therefore,
\[
\frac{A_{E,\Phi_E}(q)}
{S(W_{\bullet,\bullet}^E;q)}
=
\frac{135}{62}>1.
\]

For \(q\in E\setminus(C_2\cup L_2)\), we have \(A_{E,\Phi_E}(q)=1\) and 
\begin{align*}
    h(t)=
    \begin{dcases}
        \frac{49}{8}t^2, & t\in \left[0,\frac{2}{21}\right],\\
        \frac{1}{2}\left(\frac{3}{7}-t\right)^2, & t\in\left[\frac{2}{21},\frac{3}{7}\right].
    \end{dcases}
\end{align*}
Thus, we obtain \(S(W_{\bullet,\bullet}^E;q)=\frac{1}{9}\) and 
\[
\frac{A_{E,\Phi_E}(q)}
{S(W_{\bullet,\bullet}^E;q)}
=
9>1.
\]

The canonical divisor formula gives \(K_{S'_2}+B_{S'_2}=\psi^{\ast}K_{S_{4,2}^7}\), where \(B_{S'_2}=\frac27L_2+\frac47E\).
Applying Lemma \ref{lem:AZ-crepant-model} with the divisor
\(E\), the preceding calculations give
\[
\inf_{q\in E}
\delta_q(S'_2,B_{S'_2};\psi^{\ast}(-K_{S_{4,2}^7}))
\geq
\min\left\{
    \frac{27}{11},
    \frac{27}{10},
    \frac{135}{62},
    9
\right\}
=
\frac{135}{62}.
\]

The divisors \(L_1\subset S'_1\) and \(L_2\subset S'_2\) define the
same divisorial valuation over \(S_{4,2}^7\). Since \(\pi_2\) is an
isomorphism over \(L_2\setminus E\), pullback identifies the restricted
graded linear series and their induced filtrations at corresponding
points. Therefore,
\[
\inf_{q\in L_2\setminus E}
\delta_q(S'_2,B_{S'_2};\psi^{\ast}(-K_{S_{4,2}^7}))
\geq
\min\left\{
    \frac{75}{29},
    15,
    \frac{75}{26}
\right\}
=
\frac{75}{29}.
\]

Since \(\psi^{-1}(p)=L_2\cup E\), Lemma \ref{lem:AZ-crepant-model} yields
\[
    \delta_p(S_{4,2}^7)
    \geq
    \min\left\{
        \frac{75}{29},
        \frac{135}{62}
    \right\}
    =
    \frac{135}{62}
    >
    1.\qedhere
\]
\end{proof}

\begin{lemma}\label{lem;427 smooth}
Let \(p\in S_{4,2}^7\) be a smooth point. Then we have \(\alpha_p(S_{4,2}^7)\geq\frac{7}{10}\). In particular, \(\delta_p(S_{4,2}^7)>1\).
\end{lemma}

\begin{proof}
Apply Lemma \ref{lem:smooth-alpha-reduction} with
\((n,m,k)=(4,2,7)\). The target of the elementary transformation is
\(S_{3,2}^6\). By Lemma \ref{lem;326 smooth}, we obtain \(\alpha_q(S_{3,2}^6)\geq\frac{7}{10}\) for every smooth point \(q\in S_{3,2}^6\). Hence,
\(\alpha_p(S_{4,2}^7)\geq\frac{7}{10}\).

By Theorem \ref{thm:alpha delta}, we have
\[
    \delta_p(S_{4,2}^7)
    \geq
    \frac32\alpha_p(S_{4,2}^7)
    \geq
    \frac{21}{20}
    >
    1.\qedhere
\]
\end{proof}

\begin{theorem}\label{thm:427 stable}
    The surface \(S_{4,2}^7\) is K-stable.
\end{theorem}
\begin{proof}
Lemmas \ref{lem;427 sing} and \ref{lem;427 smooth} give \(\delta(S_{4,2}^7)\geq\frac{21}{20}>1\).
Therefore, Theorems \ref{thm:delta} and
\ref{thm:K-stable-uniform} imply that \(S_{4,2}^7\) is K-stable.
\end{proof}

Next, we deal with the surfaces \(S_{n,m}^{n+2}\). Let \(\pi_1\colon S^{(2)}_1\to \P(1,1,n)\) be the blow-up of \(\P(1,1,n)\) at \(m\) distinct smooth points \(p_1,\dots,p_m\) on a curve \(\ell\in |\mathcal O_{\P(1,1,n)}(1)|\) and at \(n+2\) general smooth  points \(q_1,\dots,q_{n+2}\). For each \(i\), there exists a unique irreducible curve \(c_i\in |\mathcal{O}_{\P(1,1,n)}(n)|\) passing through all \(q_j\) with \(j\neq i\). Moreover, we have the following intersection numbers:
    \[
    \ell^2=\frac{1}{n},\quad \ell\cdot c_i=1,\quad c_i^2=n, \quad c_i\cdot c_j=n \qquad\text{for }i\neq j.
    \]

Contracting the strict transform \(L\) of \(\ell\) on \(S^{(2)}_1\), we obtain the birational morphism \(\varphi\colon S^{(2)}_1\to S_{n,m}^{n+2}\). Let \(\pi_2\colon S^{(2)}_2\to S^{(2)}_1\) be the weighted blow-up of
\(S^{(2)}_1\) at the singular point \(o\in L\) with weights \((1,1)\). These birational morphisms are shown in the following diagram. 
\begin{figure}[H]
  \centering
  \begin{tikzpicture}
  [-,auto,node distance=1.5cm, thick,main node/.style={circle,draw,font=\sffamily \Large\bfseries}]
    \node[text=black] (1) {$\mathbb{P}(1,1,n)$};
    \node[text=black] (2) [right of=1, above of=1] {$S^{(2)}_1$};
    \node[text=black] (3) [above of=2] {$S^{(2)}_2$};
    \node[text=black] (4) [right of=2, below of=2] {$S_{n,m}^{n+2}$};
   
    \path[every node/.style={font=\sffamily}]
      (2) edge[->] node [above, left=0.5mm, pos=0.3] {\(\pi_1\)} (1)
      (3) edge[->] node [left] {\(\pi_2\)} (2)
      (2) edge[->] node [above, right=0.5mm, pos=0.3] {\(\varphi\)} (4);
  \end{tikzpicture}
    \caption{}\label{fig:general n,m}
\end{figure}
Throughout the remainder of this subsection, put \(\psi_{n,m} \coloneqq \varphi\circ\pi_2\). On the surface \(S^{(2)}_1\), let \(L\) and \(C_i\) be the strict transforms of \(\ell\) and \(c_i\), respectively. Let \(E_j\) be the exceptional curve over \(p_j\) for \(1\leq j\leq m\). Then we have the following intersection numbers:
\begin{align*}
&L^2=\frac1n-m,\qquad
L\cdot C_i=1,\qquad
C_i^2=-1,\qquad
C_i\cdot C_j=0\quad\text{for }i\neq j,\qquad
E_j^2=-1,\\
&L\cdot E_j=1,\qquad
E_i\cdot E_j=0\quad\text{for }i\neq j,\qquad
C_i\cdot E_j=0\quad\text{for all }i,j.
\end{align*}

\begin{lemma}\label{lem;L smooth}
We have
\[
\frac{A_{S_{n,m}^{n+2}}(L)}
{S_{S_{n,m}^{n+2}}(L)}
=
\frac{3(n+1)^2}{mn+2n^2+4n+1}.
\]
Moreover, we have
\[
\inf_{q\in L_{\mathrm{sm}}}
\frac{A_{L,\Phi_L}(q)}{S(W_{\bullet,\bullet}^L;q)}
\geq
\min\left\{
\frac{3(n+1)^2}{mn+2n^2+4n+1},
\frac{3n}{n+1}
\right\},
\]
where \(L_{\mathrm{sm}}\coloneqq L\setminus \Sing\left(S^{(2)}_1\right)\) and \(\Phi_L\) is the different on \(L\).
\end{lemma}
\begin{proof}
    We have
\[
-K_{S^{(2)}_1}\equiv \frac{n+2}{n+1}L+\frac{1}{n+1}(C_1+\cdots+C_{n+2})+\frac{1}{n+1}(E_1+\cdots+E_m),
\]
and
\[
K_{S^{(2)}_1}=\varphi^{\ast}(K_{S_{n,m}^{n+2}})+\left(-1+\frac{n+1}{mn-1}\right)L.
\]
Hence, \(A_{S_{n,m}^{n+2}}(L)=\frac{n+1}{mn-1}\), and we obtain 
    \[
    \varphi^{\ast}(-K_{S_{n,m}^{n+2}})-tL\equiv \left(\frac{1}{n+1}+\frac{n+1}{mn-1}-t\right)L+\frac{1}{n+1}(C_1+\cdots+C_{n+2}+E_1+\cdots+E_m).
    \]
Since the intersection matrix of \(C_1,\dots,C_{n+2},E_1,\dots,E_m\) is negative definite, the pseudoeffective threshold of \(L\) is \(\frac{1}{n+1}+\frac{n+1}{mn-1}\). Set 
\begin{equation*}
    I_1 \coloneqq \left[0,\frac{n+1}{mn-1}\right],\quad I_2 \coloneqq \left[\frac{n+1}{mn-1},\frac{1}{n+1}+\frac{n+1}{mn-1}\right].
\end{equation*}
The positive and negative parts \(P(t)\) and \(N(t)\) of the Zariski decomposition of \(\varphi^{\ast}(-K_{S_{n,m}^{n+2}})-tL\) are as follows:
    \begin{align*}
    P(t)=
    \begin{dcases}
        \varphi^{\ast}(-K_{S_{n,m}^{n+2}})-tL, & t\in I_1,\\
        \left(\frac{1}{n+1}+\frac{n+1}{mn-1}-t\right)(L+C_1+\cdots+C_{n+2}+E_1+\cdots+E_m), & t\in I_2
    \end{dcases}
    \end{align*}
    and
    \begin{align*}
        N(t)=
        \begin{dcases}
            0, & t\in I_1,\\
            \left(-\frac{n+1}{mn-1}+t\right)(C_1+\cdots+C_{n+2}+E_1+\cdots+E_m), & t\in I_2.
        \end{dcases}
    \end{align*}
    Therefore, we have
    \begin{align*}
        P(t)^2=
        \begin{dcases}
            \frac{m+n+2}{mn-1}+\left(\frac{1}{n}-m\right)t^2, &t\in I_1,\\
            \frac{(n+1)^2}{n}\left(\frac{1}{n+1}+\frac{n+1}{mn-1}-t\right)^2, & t\in I_2.
        \end{dcases}
    \end{align*}
    Hence, we obtain 
    \[
    S(L)=\frac{mn+2n^2+4n+1}{3(mn-1)(n+1)}.
    \]

    We next compute the invariant \(S(W_{\bullet,\bullet}^L;q)\) for points \(q\in L_{\mathrm{sm}}\). Set
\[
L_{CE}
\coloneqq
L\cap
\left(
\bigcup_{i=1}^{n+2}C_i
\cup
\bigcup_{j=1}^{m}E_j
\right).
\]
By the generality assumption, the points of \(L_{CE}\) are pairwise distinct, and
each of them lies on exactly one curve among the \(C_i\) and \(E_j\).

If \(q\notin L_{CE}\), then we have 
    \[
    h(t)=
    \begin{dcases}
        \frac{1}{2}\left(\frac{1}{n}-m\right)^2t^2, & t\in I_1,\\
        \frac{1}{2}\left(\frac{1}{n+1}+\frac{n+1}{mn-1}-t\right)^2\left(\frac{(n+1)^2}{n}\right)^2, & t\in I_2.
    \end{dcases}
    \]
    If \(q\in L_{CE}\), then we have
\[
h(t)=
\begin{dcases}
\frac12\left(\frac1n-m\right)^2t^2,
& t\in I_1,\\[1ex]
\begin{aligned}
&
\frac{(n+1)^2}{n}
\left(
\frac1{n+1}+\frac{n+1}{mn-1}-t
\right)
\left(
-\frac{n+1}{mn-1}+t
\right)\\
&\quad+
\frac12
\left(
\frac1{n+1}+\frac{n+1}{mn-1}-t
\right)^2
\left(
\frac{(n+1)^2}{n}
\right)^2
\end{aligned},
& t\in I_2.
\end{dcases}
\]
    
    Hence, we obtain
    \[
    S(W_{\bullet,\bullet}^L;q) =
    \begin{dcases}
         \frac{n+1}{3n}, & q\notin L_{CE},\\
         \frac{n^3+(m+4)n^2+(5+3m)n+m+1}{3n(n+1)(m+n+2)}, & q\in L_{CE}. 
    \end{dcases}
    \]

Since \(\varphi\) contracts only \(L\), none of the curves
\(C_1,\dots,C_{n+2},E_1,\dots,E_m\) contributes to the different on
\(L_{\mathrm{sm}}\). Hence, \(A_{L,\Phi_L}(q)=1\)  for every \(q\in L_{\mathrm{sm}}\).

Moreover, 
\begin{align*}
    &\frac{mn+2n^2+4n+1}{3(n+1)^2}
-
\frac{n^3+(m+4)n^2+(5+3m)n+m+1}{3n(n+1)(m+n+2)}\\
&=
\frac{
m^2n^2+2mn^3+2mn^2-3mn-m+n^4+3n^3-4n-1
}
{3n(n+1)^2(m+n+2)}
\geq 0
\end{align*}
for all \(m,n\geq2\). Indeed, if \(N_{m,n}\) denotes the numerator, then
\[
N_{m+1,n}-N_{m,n}
=
(2m+1)n^2+2n^3+2n^2-3n-1>0
\]
for all \(m,n\geq2\). Hence, it is enough to consider \(m=2\). In this
case,
\[
N_{2,n}=n^4+7n^3+8n^2-10n-3>0
\]
for every \(n\geq2\). If
\(q\notin L_{CE}\), then
\[
    \frac{A_{L,\Phi_L}(q)}
    {S(W_{\bullet,\bullet}^L;q)}
    =
    \frac{3n}{n+1}.
\]
If
\(q\in L_{CE}\), then the preceding
comparison gives
\[
    \frac{A_{L,\Phi_L}(q)}
    {S(W_{\bullet,\bullet}^L;q)}
    \geq
    \frac{3(n+1)^2}{mn+2n^2+4n+1}.
\]
Taking the infimum over \(q\in L_{\mathrm{sm}}\) gives
\[
\inf_{q\in L_{\mathrm{sm}}}
\frac{A_{L,\Phi_L}(q)}
{S(W_{\bullet,\bullet}^L;q)}
\geq
\min\left\{
    \frac{3(n+1)^2}{mn+2n^2+4n+1},
    \frac{3n}{n+1}
\right\}.\qedhere
\]
\end{proof}

\bigskip

We now compute \(\delta_p(S_{n,m}^{n+2})\) when \(m=2\), where
\(p\in S_{n,m}^{n+2}\) is the singular point. Then
\begin{align*}
K_{S^{(2)}_2}
&=
\pi_2^{\ast}(K_{S^{(2)}_1})
+\left(-1+\frac{2}{n}\right)E\\
&=
\psi_{n,m}^{\ast}(K_{S_{n,m}^{n+2}})
+\left(-1+\frac{n+1}{mn-1}\right)\widetilde L
+\left(-1+\frac{m+1}{mn-1}\right)E,
\end{align*}
where the maps \(\pi_2\) and \(\varphi\) are defined in Figure \ref{fig:general n,m}.
Therefore, we have \(A_{S_{n,m}^{n+2}}(E)=\frac{m+1}{mn-1}\).

When \(m=2\), we only need to consider three cases: \((n,m,k)=(3,2,5),(4,2,6),(5,2,7)\). The surface \(S_{n,2}^{n+2}\) is obtained by blowing up
two distinct smooth points \(p_1,p_2\) on a member
\(\ell\in|\mathcal O_{\P(1,1,n)}(1)|\) and \(n+2\) general smooth
points \(q_1,\dots,q_{n+2}\) outside \(\ell\), and then contracting
the strict transform of \(\ell\). Let
\(\mu_n\colon\mathbb F_n\to\P(1,1,n)\) be the minimal resolution, let \(B_n\) be its negative section, and let \(\mathfrak f\) be the fibre class. The strict transform of \(\ell\) is a fibre, which we denote by \(\mathfrak f_0\).

For \(j\in\{1,2\}\), let \(D_j\in|B_n+(n+1)\mathfrak f|\) be the unique member passing through \(p_j\) and \(q_1,\dots,q_{n+2}\). The divisor \(B_n+(n+1)\mathfrak f\) is very ample, and its complete linear system has projective dimension \(n+3\). The chosen \(n+3\) general points impose independent conditions.

Let \(F_i\) be the fibre passing through \(q_i\), and let
\(\Gamma_{n,2}\) be the image on \(\P(1,1,n)\) of the reduced divisor \(D_1+D_2+F_1+\cdots+F_{n+2}\).
Then we have \(\Gamma_{n,2} \in |\mathcal O_{\P(1,1,n)}(3n+4)|\).
It has multiplicity \(3\) at every \(q_i\) and multiplicity \(1\) at
each \(p_j\). For a general configuration, the strict transforms of
the irreducible components are pairwise disjoint \((-1)\)-curves. Indeed, after imposing the \(n+2\) conditions at \(q_1,\dots,q_{n+2}\), the linear system \(|B_n+(n+1)\mathfrak f|\) becomes a pencil. Two general members of this pencil are smooth and
irreducible, intersect transversely precisely at
\(q_1,\dots,q_{n+2}\), and meet the fibre \(\mathfrak f_0\) at two
distinct points. Therefore, the required transversality and
disjointness properties hold on a nonempty open subset of the
configuration space.

Then
\[
-K_{S_2^{(2)}}
\sim_{\Q}
\frac43E+\frac23\widetilde L
+\frac13\widetilde\Gamma_{n,2},
\]
where \(E\) is the exceptional curve, while
\(\widetilde L\) and \(\widetilde\Gamma_{n,2}\) are the strict
transforms of \(\ell\) and \(\Gamma_{n,2}\), respectively. Moreover, we have the following intersection numbers:
\[
E^2=-n,\qquad
\widetilde L^2=-2,\qquad
\widetilde\Gamma_{n,2}^2=-(n+4),\qquad
E\cdot\widetilde\Gamma_{n,2}=n+4,\qquad
E\cdot\widetilde L=1,\qquad
\widetilde L\cdot\widetilde\Gamma_{n,2}=0.
\]
Here, \(\widetilde\Gamma_{n,2}\) denotes the sum of its pairwise
disjoint irreducible components.

Note that we have
\begin{equation}\label{eq:m2-pullback}
\psi_{n,2}^{\ast}(-K_{S_{n,2}^{n+2}})-tE
=
\left(\frac13+\frac3{2n-1}-t\right)E
+
\left(-\frac13+\frac{n+1}{2n-1}\right)\widetilde L
+
\frac13\widetilde\Gamma_{n,2}.
\end{equation}

Since the intersection matrix of \(\widetilde L\) and \(\widetilde\Gamma_{n,2}\) is negative definite, we have \(\tau(E)=\frac{1}{3}+\frac{3}{2n-1}\).

Now, we estimate the \(\delta\)-invariant case by case.

\subsubsection{Case: \(S_{3,2}^5\)}

\begin{lemma}\label{lem;325 sing}
    Let \(p\in S_{3,2}^5\) be the singular point of type \(\frac{1}{5}(1,3)\). Then we have \(\delta_p(S_{3,2}^5)>1\).
\end{lemma}
\begin{proof}
    We note that \(A_{S_{3,2}^5}(E)=\frac{3}{5}\) and by \eqref{eq:m2-pullback}, we have
\begin{align*}
    \psi_{3,2}^{\ast}(-K_{S_{3,2}^5})-tE=\left(\frac{14}{15}-t\right)E+\frac{7}{15}\widetilde L+\frac{1}{3}\widetilde\Gamma_{3,2}.
\end{align*}
The positive and negative parts of the Zariski decomposition of \(\psi_{3,2}^{\ast}(-K_{S_{3,2}^5})-tE\) are as follows:
\begin{align*}
    P(t)=
    \begin{dcases}
        \left(\frac{14}{15}-t\right)\left(E+\frac{1}{2}\widetilde L\right)+\frac{1}{3}\widetilde\Gamma_{3,2}, & t\in\left[0,\frac{3}{5}\right],\\
        \left(\frac{14}{15}-t\right)\left(E+\frac{1}{2}\widetilde L+\widetilde\Gamma_{3,2}\right), & t\in\left[\frac{3}{5},\frac{14}{15}\right],
    \end{dcases}
\end{align*}
and
\begin{align*}
    N(t)=
    \begin{dcases}
        \frac{1}{2}t\widetilde L, & t\in\left[0,\frac{3}{5}\right],\\
        \frac{1}{2}t\widetilde L+\left(t-\frac{3}{5}\right)\widetilde\Gamma_{3,2}, & t\in\left[\frac{3}{5},\frac{14}{15}\right].
    \end{dcases}
\end{align*}

Therefore, we obtain \(S_{S_{3,2}^5}(E)=\frac{5}{7}\left(\frac{33}{50}+\frac{1}{18}\right)=\frac{23}{45}\).
Hence, 
\[
\frac{A_{S_{3,2}^5}(E)}{S_{S_{3,2}^5}(E)}=\dfrac{27}{23}.
\]

In order to apply the Abban--Zhuang theory, we need to compute the following intersection number
\begin{align*}
    P(t)\cdot E=
    \begin{dcases}
        \frac{5}{2}t, & t\in\left[0,\frac{3}{5}\right],\\
        \frac{9}{2}\left(\frac{14}{15}-t\right), & t\in\left[\frac{3}{5},\frac{14}{15}\right].
    \end{dcases}
\end{align*}

We note that \(\widetilde L \cap \Supp(\widetilde\Gamma_{3,2}) = \emptyset\). For \(q\in E\setminus \left(\Supp(\widetilde\Gamma_{3,2})\cup\widetilde L \right)\), we have
\begin{align*}
    h(t)=
    \begin{dcases}
        \frac{25}{8}t^2, & t\in\left[0,\frac{3}{5}\right],\\
        \frac{81}{8}\left(\frac{14}{15}-t\right)^2, & t\in\left[\frac{3}{5},\frac{14}{15}\right].
    \end{dcases}
\end{align*}
Hence, \(S(W_{\bullet,\bullet}^E;q) = \frac{10}{7}\left(\frac{9}{40}+\frac{1}{8}\right) = \frac{1}{2}\).

For \(q\in E\cap \widetilde L\), we have
\begin{align*}
    h(t)=
    \begin{dcases}
        \frac{5}{4}t^2+\frac{25}{8}t^2, & t\in\left[0,\frac{3}{5}\right],\\
        \frac{9}{4}t\left(\frac{14}{15}-t\right)+\frac{81}{8}\left(\frac{14}{15}-t\right)^2, & t\in\left[\frac{3}{5},\frac{14}{15}\right].
    \end{dcases}
\end{align*}
Hence, \(S(W_{\bullet,\bullet}^E;q) = \frac{10}{7}\left(\frac{63}{200}+\frac{77}{360}\right) = \frac{34}{45}\). At \(q=E\cap\widetilde L\), we have \(A_{E,\Phi_E}(q)=A_{S_{3,2}^5}(\widetilde L)=\frac{4}{5}\). Thus, we have
\[
\frac{A_{E,\Phi_E}(q)}
{S(W_{\bullet,\bullet}^E;q)}
=
\frac{18}{17}>1.
\]

For \(q\in E\cap \Supp(\widetilde\Gamma_{3,2})\), since \((E\cdot \widetilde\Gamma_{3,2})_q\leq E\cdot \widetilde\Gamma_{3,2}=7\), we have
\[
h(t)=\dfrac{25}{8}t^2, \quad t\in\left[0,\frac{3}{5}\right],
\]
and
\[
h(t)\leq \dfrac{63}{2}\left(\dfrac{14}{15}-t\right)\left(t-\dfrac{3}{5}\right)+\dfrac{81}{8}\left(\dfrac{14}{15}-t\right)^2, \quad t\in\left[\frac{3}{5},\frac{14}{15}\right].
\]
Hence, \(S(W_{\bullet,\bullet}^E;q) \leq \frac{10}{7}\left(\frac{9}{40}+\frac{23}{72}\right) = \frac{7}{9}\).
Since the irreducible component of \(\widetilde\Gamma_{3,2}\) passing
through \(q\) is not exceptional over \(S_{3,2}^5\), no component of
the different contributes at \(q\). Hence, we have \(A_{E,\Phi_E}(q)=1\) and thus, 
\[
\frac{A_{E,\Phi_E}(q)}
{S(W_{\bullet,\bullet}^E;q)}
\geq
\frac{9}{7}>1.
\]

The canonical divisor formula gives \(K_{S_2^{(2)}}+B_{S_2^{(2)}}=\psi_{3,2}^{\ast}K_{S_{3,2}^5}\), where \(B_{S_2^{(2)}}=\frac15\widetilde L+\frac25E\).
Applying Lemma \ref{lem:AZ-crepant-model} with the divisor
\(E\), the preceding calculations give
\[
\inf_{q\in E}
\delta_q(S_2^{(2)},B_{S_2^{(2)}};\psi_{3,2}^{\ast}(-K_{S_{3,2}^5}))
\geq
\min\left\{
    \frac{27}{23},
    2,
    \frac{18}{17},
    \frac97
\right\}
=
\frac{18}{17}.
\]

The morphism \(\pi_2\) is an isomorphism over
\(\widetilde L\setminus E\), and identifies this open subset with
\(L_{\mathrm{sm}}\subset S_1^{(2)}\). Under this identification, the
restricted graded linear series and their induced filtrations agree.
Hence, Lemmas \ref{lem;L smooth} and
\ref{lem:AZ-crepant-model} give
\[
\inf_{q\in\widetilde L\setminus E}
\delta_q(S_2^{(2)},B_{S_2^{(2)}};\psi_{3,2}^{\ast}(-K_{S_{3,2}^5}))
\geq
\min\left\{
    \frac{48}{37},
    \frac94
\right\}
=
\frac{48}{37}.
\]
Since \(\psi_{3,2}^{-1}(p)=\widetilde L\cup E\), Lemma \ref{lem:AZ-crepant-model} yields
\[
    \delta_p(S_{3,2}^5)
    \geq
    \min\left\{
        \frac{48}{37},
        \frac{18}{17}
    \right\}
    =
    \frac{18}{17}
    >
    1.\qedhere
\]
\end{proof}

\begin{lemma}\label{lem;325 smooth}
Let \(p\in S_{3,2}^5\) be a smooth point. Then we have \(\alpha_p(S_{3,2}^5)\geq\frac{7}{10}\). In particular, \(\delta_p(S_{3,2}^5)>1\).
\end{lemma}

\begin{proof}
Apply Lemma \ref{lem:smooth-alpha-reduction} with
\((n,m,k)=(3,2,5)\). The target of the elementary transformation is
\(S_{2,2}^4\). By Lemma \ref{lem:base-22-smooth-alpha}, we obtain \(\alpha_q(S_{2,2}^4)\geq\frac34>\frac{7}{10}\) for every smooth point \(q\in S_{2,2}^4\). Hence,
\(\alpha_p(S_{3,2}^5)\geq\frac{7}{10}\).

By Theorem \ref{thm:alpha delta}, we have
\[
    \delta_p(S_{3,2}^5)
    \geq
    \frac32\alpha_p(S_{3,2}^5)
    \geq
    \frac{21}{20}
    >
    1.\qedhere
\]
\end{proof}

\begin{theorem}\label{thm:325 stable}
    The surface \(S_{3,2}^5\) is K-stable.
\end{theorem}
\begin{proof}
Lemmas \ref{lem;325 sing} and \ref{lem;325 smooth} give \(\delta(S_{3,2}^5)\geq\frac{21}{20}>1\). Therefore, Theorems \ref{thm:delta} and \ref{thm:K-stable-uniform} imply that \(S_{3,2}^5\) is K-stable.
\end{proof}

\subsubsection{Case: \(S_{4,2}^6\)}

\begin{lemma}\label{lem;426 sing}
    Let \(p\in S_{4,2}^6\) be the singular point of type \(\frac{1}{7}(1,4)\). Then we have \(\delta_p(S_{4,2}^6)>1\).
\end{lemma}
\begin{proof}
    We note that \(A_{S_{4,2}^6}(E)=\frac{3}{7}\) and by \eqref{eq:m2-pullback}, we have
\begin{align*}
    \psi_{4,2}^{\ast}(-K_{S_{4,2}^6})-tE=\left(\frac{16}{21}-t\right)E+\frac{8}{21}\widetilde L+\frac{1}{3}\widetilde\Gamma_{4,2}.
\end{align*}

The positive and negative parts of the Zariski decomposition are as follows:
\begin{align*}
    P(t)=
    \begin{dcases}
        \left(\frac{16}{21}-t\right)\left(E+\frac{1}{2}\widetilde L\right)+\frac{1}{3}\widetilde\Gamma_{4,2}, & t\in\left[0,\frac{3}{7}\right],\\
        \left(\frac{16}{21}-t\right)\left(E+\frac{1}{2}\widetilde L+\widetilde\Gamma_{4,2}\right), & t\in\left[\frac{3}{7},\frac{16}{21}\right],
    \end{dcases}
\end{align*}
and
\begin{align*}
    N(t)=
    \begin{dcases}
        \frac{1}{2}t\widetilde L, & t\in\left[0,\frac{3}{7}\right],\\
        \frac{1}{2}t\widetilde L+\left(t-\frac{3}{7}\right)\widetilde\Gamma_{4,2}, & t\in\left[\frac{3}{7},\frac{16}{21}\right].
    \end{dcases}
\end{align*}

Therefore, we obtain \(S_{S_{4,2}^6}(E)=\frac{7}{8}\left(\frac{39}{98}+\frac{1}{18}\right)=\frac{25}{63}\). Hence,
\[
\frac{A_{S_{4,2}^6}(E)}{S_{S_{4,2}^6}(E)}=\frac{27}{25}.
\]
Moreover, we have
\begin{align*}
    P(t)\cdot E=
    \begin{dcases}
        \frac{7}{2}t, & t\in\left[0,\frac{3}{7}\right],\\
        \frac{9}{2}\left(\frac{16}{21}-t\right), & t\in\left[\frac{3}{7},\frac{16}{21}\right].
    \end{dcases}
\end{align*}

For \(q\in E\setminus \left(\Supp(\widetilde\Gamma_{4,2})\cup\widetilde L \right)\), we have
\begin{align*}
    h(t)=
    \begin{dcases}
        \frac{49}{8}t^2, & t\in\left[0,\frac{3}{7}\right],\\
        \frac{81}{8}\left(\frac{16}{21}-t\right)^2, & t\in\left[\frac{3}{7},\frac{16}{21}\right].
    \end{dcases}
\end{align*}
Hence, \(S(W_{\bullet,\bullet}^E;q) = \frac{7}{4}\left(\frac{9}{56}+\frac{1}{8}\right) = \frac{1}{2}\).
For \(q\in E\cap \widetilde L\), we have
\begin{align*}
    h(t)=
    \begin{dcases}
        \frac{7}{4}t^2+\frac{49}{8}t^2, & t\in\left[0,\frac{3}{7}\right],\\
        \frac{9}{4}t\left(\frac{16}{21}-t\right)+\frac{81}{8}\left(\frac{16}{21}-t\right)^2, & t\in\left[\frac{3}{7},\frac{16}{21}\right].
    \end{dcases}
\end{align*}
Hence, \(S(W_{\bullet,\bullet}^E;q) = \frac{7}{4}\left(\frac{81}{392}+\frac{97}{504}\right) = \frac{44}{63}\). At \(q=E\cap\widetilde L\), we have \(A_{E,\Phi_E}(q)=A_{S_{4,2}^6}(\widetilde L)=\frac{5}{7}\).
Thus, we have
\[
\frac{A_{E,\Phi_E}(q)}
{S(W_{\bullet,\bullet}^E;q)}
=
\frac{45}{44}>1.
\]

For \(q\in E\cap \Supp(\widetilde\Gamma_{4,2})\), since \((E\cdot \widetilde\Gamma_{4,2})_q\leq E\cdot \widetilde\Gamma_{4,2}=8\), we have
\[
h(t)=\dfrac{49}{8}t^2, \quad t\in\left[0,\frac{3}{7}\right],
\]
and
\[
h(t)\leq \dfrac{72}{2}\left(\dfrac{16}{21}-t\right)\left(t-\dfrac{3}{7}\right)+\dfrac{81}{8}\left(\dfrac{16}{21}-t\right)^2, \quad t\in\left[\frac{3}{7},\frac{16}{21}\right].
\]
Hence, \(S(W_{\bullet,\bullet}^E;q) \leq \frac{7}{4}\left(\frac{9}{56}+\frac{25}{72}\right) = \frac{8}{9}\).
Since the irreducible component of \(\widetilde\Gamma_{4,2}\) passing
through \(q\) is not exceptional over \(S_{4,2}^6\), we have \(A_{E,\Phi_E}(q)=1\) and thus,
\[
\frac{A_{E,\Phi_E}(q)}
{S(W_{\bullet,\bullet}^E;q)}
\geq
\frac{9}{8}>1.
\]

The canonical divisor formula gives \(K_{S_2^{(2)}}+B_{S_2^{(2)}}=\psi_{4,2}^{\ast}K_{S_{4,2}^6}\), where \(B_{S_2^{(2)}}=\frac27\widetilde L+\frac47E\).
Applying Lemma \ref{lem:AZ-crepant-model} with the divisor
\(E\), the preceding calculations give
\[
\inf_{q\in E}
\delta_q(S_2^{(2)},B_{S_2^{(2)}};\psi_{4,2}^{\ast}(-K_{S_{4,2}^6}))
\geq
\min\left\{
    \frac{27}{25},
    2,
    \frac{45}{44},
    \frac98
\right\}
=
\frac{45}{44}.
\]

The morphism \(\pi_2\) is an isomorphism over
\(\widetilde L\setminus E\), and identifies this open subset with
\(L_{\mathrm{sm}}\subset S_1^{(2)}\). Under this identification, the
restricted graded linear series and their induced filtrations agree.
Hence, Lemmas \ref{lem;L smooth} and
\ref{lem:AZ-crepant-model} give
\[
\inf_{q\in\widetilde L\setminus E}
\delta_q(S_2^{(2)},B_{S_2^{(2)}};\psi_{4,2}^{\ast}(-K_{S_{4,2}^6}))
\geq
\min\left\{
    \frac{25}{19},
    \frac{12}{5}
\right\}
=
\frac{25}{19}.
\]
Since \(\psi_{4,2}^{-1}(p)=\widetilde L\cup E\), Lemma \ref{lem:AZ-crepant-model} yields
\[
    \delta_p(S_{4,2}^6)
    \geq
    \min\left\{
        \frac{25}{19},
        \frac{45}{44}
    \right\}
    =
    \frac{45}{44}
    >
    1.\qedhere
\]
\end{proof}

Now we compute \(\delta_p(S_{4,2}^6)\) for a smooth point \(p\).

\begin{lemma}\label{lem;426 smooth}
Let \(p\in S_{4,2}^6\) be a smooth point. Then we have \(\alpha_p(S_{4,2}^6)\geq\frac{7}{10}\). In particular, \(\delta_p(S_{4,2}^6)>1\).
\end{lemma}

\begin{proof}
Apply Lemma \ref{lem:smooth-alpha-reduction} with
\((n,m,k)=(4,2,6)\). The target of the elementary transformation is
\(S_{3,2}^5\). By Lemma \ref{lem;325 smooth}, we obtain \(\alpha_q(S_{3,2}^5)\geq\frac{7}{10}\) for every smooth point \(q\in S_{3,2}^5\). Hence,
\(\alpha_p(S_{4,2}^6)\geq\frac{7}{10}\).

By Theorem \ref{thm:alpha delta}, we have
\[
    \delta_p(S_{4,2}^6)
    \geq
    \frac32\alpha_p(S_{4,2}^6)
    \geq
    \frac{21}{20}
    >
    1.\qedhere
\]
\end{proof}

\begin{theorem}\label{thm:426 stable}
    The surface \(S_{4,2}^6\) is K-stable.
\end{theorem}
\begin{proof}
Lemmas \ref{lem;426 sing} and \ref{lem;426 smooth} give \(\delta(S_{4,2}^6)\geq\frac{45}{44}>1\).
Therefore, Theorems \ref{thm:delta} and
\ref{thm:K-stable-uniform} imply that \(S_{4,2}^6\) is K-stable.
\end{proof}

\subsubsection{Case: \(S_{5,2}^7\)}

\begin{lemma}\label{lem;527 sing}
    Let \(p\in S_{5,2}^7\) be the singular point of type \(\frac{1}{9}(1,5)\). Then we have \(\delta_p(S_{5,2}^7)=1\).
\end{lemma}
\begin{proof}
    We note that \(A_{S_{5,2}^7}(E)=\frac{1}{3}\) and by \eqref{eq:m2-pullback}, we have
\begin{align*}
    \psi_{5,2}^{\ast}(-K_{S_{5,2}^7})-tE=\left(\frac{2}{3}-t\right)E+\frac{1}{3}\widetilde L+\frac{1}{3}\widetilde\Gamma_{5,2}.
\end{align*}

The positive and negative parts of the Zariski decomposition are as follows:

\begin{align*}
    P(t)=
    \begin{dcases}
        \left(\frac{2}{3}-t\right)\left(E+\frac{1}{2}\widetilde L\right)+\frac{1}{3}\widetilde\Gamma_{5,2}, & t\in\left[0,\frac{1}{3}\right],\\
        \left(\frac{2}{3}-t\right)\left(E+\frac{1}{2}\widetilde L+\widetilde\Gamma_{5,2}\right), & t\in\left[\frac{1}{3},\frac{2}{3}\right],
    \end{dcases}
\end{align*}
and
\begin{align*}
    N(t)=
    \begin{dcases}
        \frac{1}{2}t\widetilde L, & t\in\left[0,\frac{1}{3}\right],\\
        \frac{1}{2}t\widetilde L+\left(t-\frac{1}{3}\right)\widetilde\Gamma_{5,2}, & t\in\left[\frac{1}{3},\frac{2}{3}\right].
    \end{dcases}
\end{align*}

Therefore, we obtain \(S_{S_{5,2}^7}(E)=\frac{5}{18}+\frac{1}{18}=\frac{1}{3}\). Hence,
\[
\frac{A_{S_{5,2}^7}(E)}{S_{S_{5,2}^7}(E)}=1.
\]
Moreover, we have
\begin{align*}
    P(t)\cdot E=
    \begin{dcases}
        \frac{9}{2}t, & t\in\left[0,\frac{1}{3}\right],\\
        \frac{9}{2}\left(\frac{2}{3}-t\right), & t\in\left[\frac{1}{3},\frac{2}{3}\right].
    \end{dcases}
\end{align*}

For \(q\in E\setminus \left(\Supp(\widetilde\Gamma_{5,2})\cup\widetilde L\right)\), we have
\begin{align*}
    h(t)=
    \begin{dcases}
        \frac{81}{8}t^2, & t\in\left[0,\frac{1}{3}\right],\\
        \frac{81}{8}\left(\frac{2}{3}-t\right)^2, & t\in\left[\frac{1}{3},\frac{2}{3}\right].
    \end{dcases}
\end{align*}
Hence, \(S(W_{\bullet,\bullet}^E;q) = 2\left(\frac{1}{8}+\frac{1}{8}\right) = \frac{1}{2}\).

For \(q\in E\cap \widetilde L\), we have
\begin{align*}
    h(t)=
    \begin{dcases}
        \frac{9}{4}t^2+\frac{81}{8}t^2, & t\in\left[0,\frac{1}{3}\right],\\
        \frac{9}{4}t\left(\frac{2}{3}-t\right)+\frac{81}{8}\left(\frac{2}{3}-t\right)^2, & t\in\left[\frac{1}{3},\frac{2}{3}\right].
    \end{dcases}
\end{align*}
Hence, \(S(W_{\bullet,\bullet}^E;q) = 2\left(\frac{11}{72}+\frac{13}{72}\right) = \frac{2}{3}\). At \(q=E\cap\widetilde L\), we have \(A_{E,\Phi_E}(q)=A_{S_{5,2}^7}(\widetilde L)=\frac{2}{3}\). Thus, we get
\[
\frac{A_{E,\Phi_E}(q)}
{S(W_{\bullet,\bullet}^E;q)}
=
1.
\]

For \(q\in E\cap \Supp(\widetilde\Gamma_{5,2})\), since \((E\cdot \widetilde\Gamma_{5,2})_q\leq E\cdot \widetilde\Gamma_{5,2}=9\), we have
\[
h(t)=\dfrac{81}{8}t^2, \quad t\in\left[0,\frac{1}{3}\right],
\]
and
\[
h(t)\leq \dfrac{81}{2}\left(\dfrac{2}{3}-t\right)\left(t-\dfrac{1}{3}\right)+\dfrac{81}{8}\left(\dfrac{2}{3}-t\right)^2, \quad t\in\left[\frac{1}{3},\frac{2}{3}\right].
\]
Hence, \(S(W_{\bullet,\bullet}^E;q) \leq 2\left(\frac{1}{8}+\frac{3}{8}\right) = 1\). Since the irreducible component of \(\widetilde\Gamma_{5,2}\) passing through \(q\) is not exceptional over \(S_{5,2}^7\), we have \(A_{E,\Phi_E}(q)=1\) and thus,
\[
\frac{A_{E,\Phi_E}(q)}
{S(W_{\bullet,\bullet}^E;q)}
\geq
1.
\]

The canonical divisor formula gives \(K_{S_2^{(2)}}+B_{S_2^{(2)}}=\psi_{5,2}^{\ast}K_{S_{5,2}^7}\), where \(B_{S_2^{(2)}}=\frac13\widetilde L+\frac23E\).
Applying Lemma \ref{lem:AZ-crepant-model} with the divisor
\(E\), the preceding calculations give
\[
\inf_{q\in E}
\delta_q(S_2^{(2)},B_{S_2^{(2)}};\psi_{5,2}^{\ast}(-K_{S_{5,2}^7}))
\geq
\min\left\{
    1,
    2,
    1,
    1
\right\}
=
1.
\]

The morphism \(\pi_2\) is an isomorphism over
\(\widetilde L\setminus E\), and identifies this open subset with
\(L_{\mathrm{sm}}\subset S_1^{(2)}\). Under this identification, the
restricted graded linear series and their induced filtrations agree.
Hence, Lemmas \ref{lem;L smooth} and
\ref{lem:AZ-crepant-model} give
\[
\inf_{q\in\widetilde L\setminus E}
\delta_q(S_2^{(2)},B_{S_2^{(2)}};\psi_{5,2}^{\ast}(-K_{S_{5,2}^7}))
\geq
\min\left\{
    \frac43,
    \frac52
\right\}
=
\frac43.
\]

Since \(\psi_{5,2}^{-1}(p)=\widetilde L\cup E\), Lemma \ref{lem:AZ-crepant-model} yields \(\delta_p(S_{5,2}^7)\geq1\). On the other hand, the divisor \(E\) gives
\[
    \delta_p(S_{5,2}^7)
    \leq
    \frac{A_{S_{5,2}^7}(E)}
    {S_{S_{5,2}^7}(E)}
    =
    1.
\]
Consequently, \(\delta_p(S_{5,2}^7)=1\).
\end{proof}

\begin{lemma}\label{lem;527 smooth}
Let \(p\in S_{5,2}^7\) be a smooth point. Then we have \(\alpha_p(S_{5,2}^7)\geq\frac{7}{10}\). In particular, \(\delta_p(S_{5,2}^7)>1\).
\end{lemma}

\begin{proof}
Apply Lemma \ref{lem:smooth-alpha-reduction} with
\((n,m,k)=(5,2,7)\). The target of the elementary transformation is
\(S_{4,2}^6\). By Lemma \ref{lem;426 smooth}, we obtain \(\alpha_q(S_{4,2}^6)\geq\frac{7}{10}\) for every smooth point \(q\in S_{4,2}^6\). Hence,
\(\alpha_p(S_{5,2}^7)\geq\frac{7}{10}\).

By Theorem \ref{thm:alpha delta}, we have
\[
    \delta_p(S_{5,2}^7)
    \geq
    \frac32\alpha_p(S_{5,2}^7)
    \geq
    \frac{21}{20}
    >
    1.\qedhere
\]
\end{proof}

We claim that the surface \(S_{5,2}^7\) is strictly K-semistable by showing that the automorphism group of \(S_{5,2}^7\) is finite.

\begin{lemma}\label{lem:335-527-aut-fin}
The automorphism groups \(\Aut(S_{3,3}^5)\) and \(\Aut(S_{5,2}^7)\) are finite.
\end{lemma}

\begin{proof}
Let \((n,m)\in \{(3,3), (5,2)\}\), and let \(X\coloneqq S_{n,m}^{n+2}\). It suffices to prove that \(\Aut^0(X)\) is trivial.

Let \(\nu\colon\widetilde X\to X\) be the minimal resolution. Every automorphism of \(X\) lifts to \(\widetilde X\). Since \(\Aut^0(X)\) is connected, its action on the discrete group \(\operatorname{NS}(\widetilde X)\) is trivial.

Let \(C\subset\widetilde X\) be an irreducible curve with
\(C^2<0\). Then \(C\) is the unique irreducible curve in its numerical
class. Indeed, if \(C'\neq C\) were another irreducible curve
numerically equivalent to \(C\), then \(C\cdot C'=C^2<0\),
which is impossible for two distinct irreducible curves on a smooth
surface. Therefore, every negative curve occurring in the construction
of \(\widetilde X\) is preserved by \(\Aut^0(X)\).

Using Notation \ref{not:standard-partial-resolution}, we equivariantly contract the exceptional curves over \(a_1,\dots,a_m,b_1,\dots,b_{n+2}\).
The resulting surface is \(\mathbb F_n\). Contracting its negative
section gives an induced action of \(\Aut^0(X)\) on
\(\P(1,1,n)\). This action preserves \(\ell\) and fixes each point \(b_1,\dots,b_{n+2}\).

Every automorphism of \(\P(1,1,n)_{x,y,z}\) is induced by a graded
automorphism of the form
\[
[x:y:z]
\longmapsto
[ax+by:cx+dy:\lambda z+g_n(x,y)],
\]
where
\[
\begin{pmatrix}
a&b\\
c&d
\end{pmatrix}
\in\operatorname{GL}_2(\C),
\qquad
\lambda\in\C^{\ast},
\]
and \(g_n(x,y)\) is a homogeneous polynomial of degree \(n\).

The projections of at least three of the points \(b_i\) to
\(\P^1_{x,y}\) are distinct. Hence, an automorphism fixing every
\(b_i\) induces the identity on \(\P^1_{x,y}\). After rescaling the weighted homogeneous coordinates, it has the form \([x:y:z]\longmapsto[x:y:\lambda z+g_n(x,y)]\).

Write \(b_i=[x_i:y_i:z_i]\).
Fixing \(b_i\) imposes the condition \(g_n(x_i,y_i)+(\lambda-1)z_i=0\)
for every \(i\). These equations are given by the evaluation map
\begin{align*}
H^0\bigl(\P^1,\mathcal O_{\P^1}(n)\bigr)\oplus\C
&\longrightarrow
\C^{n+2},\\
(g,u)
&\longmapsto
\bigl(g(x_i,y_i)+uz_i\bigr)_{i=1}^{n+2}.
\end{align*}
For general points \(b_1,\dots,b_{n+2}\), this map is an
isomorphism. Indeed, its determinant is nonzero for a configuration
consisting of \(n+1\) points with \(z_i=0\) and distinct projections
to \(\P^1\), together with one point having nonzero \(z\)-coordinate.
Thus, the determinant is nonzero on a nonempty open subset of the
configuration space.

It follows that \(g_n=0\) and \(\lambda=1\).
An automorphism of \(\widetilde X\) inducing the identity on
\(\P(1,1,n)\) is the identity on a dense open subset of
\(\widetilde X\), and hence is the identity everywhere. Therefore,
\(\Aut^0(X)\) is trivial.

For a sufficiently divisible positive integer \(r\), the divisor
\(-rK_X\) is very ample. Hence, \(\Aut(X)\) is a closed subgroup of \(\operatorname{PGL}\bigl(H^0(X,-rK_X)\bigr)\),
and is an algebraic group of finite type. Since its identity component
is trivial, \(\Aut(X)\) is finite.
\end{proof}

\begin{theorem}\label{thm:527-semistable}
The surface \(S_{5,2}^7\) is strictly K-semistable.
\end{theorem}

\begin{proof}
Lemmas \ref{lem;527 sing} and \ref{lem;527 smooth} give \(\delta(S_{5,2}^7)=1\). Hence, \(S_{5,2}^7\) is K-semistable by Theorem \ref{thm:delta}, and it is not K-stable by Theorem \ref{thm:K-stable-uniform}. Its automorphism group is finite by Lemma \ref{lem:335-527-aut-fin}. Consequently, Theorem \ref{thm:finite aut} shows that it is not K-polystable. Therefore, \(S_{5,2}^7\) is strictly K-semistable.
\end{proof}

When \(m=3\), we only have two possibilities: \((n,m)=(3,3),(4,3)\). Let
\(\mu_n\colon\mathbb F_n\to\P(1,1,n)\)
be the minimal resolution, let \(B_n\) be its negative section, and let
\(\mathfrak f\) be the fibre class. The strict transform of \(\ell\) is a fibre \(\mathfrak f_0\).

For \(j\in\{1,2,3\}\), let \(D_j\in|B_n+(n+1)\mathfrak f|\) be the unique member passing through \(p_j\) and
\(q_1,\dots,q_{n+2}\). Let \(F_i\) be the fibre through \(q_i\), and
let \(\Gamma_{n,3}\) be the image on \(\P(1,1,n)\) of the reduced
divisor \(D_1+D_2+D_3+F_1+\cdots+F_{n+2}\).
Then we have \(\Gamma_{n,3} \in |\mathcal O_{\P(1,1,n)}(4n+5)|\).
It has multiplicity \(4\) at every \(q_i\) and multiplicity \(1\) at
each \(p_j\). For a general configuration, the strict transforms of
its irreducible components are pairwise disjoint \((-1)\)-curves. Indeed, after imposing the \(n+2\) conditions at \(q_1,\dots,q_{n+2}\), the linear system \(|B_n+(n+1)\mathfrak f|\) becomes a pencil. Three general members of this pencil are smooth and
irreducible, have the prescribed common base points, and meet
\(\mathfrak f_0\) at three distinct points. Therefore, the required
transversality and disjointness properties hold for a general
configuration.

Then 
\[
-K_{S_2^{(2)}}
\sim_{\Q}
\frac54E+\frac34\widetilde L
+\frac14\widetilde\Gamma_{n,3}.
\]
Moreover, we have the following intersection numbers:
\[
E^2=-n,\qquad
\widetilde L^2=-3,\qquad
\widetilde\Gamma_{n,3}^2=-(n+5),\qquad
E\cdot\widetilde\Gamma_{n,3}=n+5,\qquad
E\cdot\widetilde L=1,\qquad
\widetilde L\cdot\widetilde\Gamma_{n,3}=0.
\]
Here, \(\widetilde\Gamma_{n,3}\) denotes the sum of its pairwise
disjoint irreducible components.

Note that we have
\begin{equation}\label{eq:m3-pullback}
\psi_{n,3}^{\ast}(-K_{S_{n,3}^{n+2}})-tE
=
\left(\frac14+\frac4{3n-1}-t\right)E
+
\left(-\frac14+\frac{n+1}{3n-1}\right)\widetilde L
+
\frac14\widetilde\Gamma_{n,3},
\end{equation}
where the maps are defined in Figure \ref{fig:general n,m}.
Since the intersection matrix of \(\widetilde L\) and \(\widetilde\Gamma_{n,3}\) is negative definite, we have \(\tau(E)=\frac{1}{4}+\frac{4}{3n-1}\).

\subsubsection{Case: \(S_{3,3}^5\)}

\begin{lemma}\label{lem;335 sing}
Let \(p\in S_{3,3}^5\) be the singular point of type
\(\frac18(1,3)\). Then we have \(\delta_p(S_{3,3}^5)=1\).
\end{lemma}

\begin{proof}
The canonical divisor formula gives \(K_{S_2^{(2)}}+B_{S_2^{(2)}}=\psi_{3,3}^{\ast}K_{S_{3,3}^5}\), where \(B_{S_2^{(2)}}=\frac12\widetilde L+\frac12E\).

By Lemmas \ref{lem;L smooth} and
\ref{lem:AZ-crepant-model}, we have
\[
\inf_{z\in\widetilde L\setminus E}
\delta_z(S_2^{(2)},B_{S_2^{(2)}};\psi_{3,3}^{\ast}(-K_{S_{3,3}^5}))
\geq
\min\left\{
    \frac65,
    \frac94
\right\}
=
\frac65.
\]

By \eqref{eq:m3-pullback}, we have
\[
\psi_{3,3}^{\ast}(-K_{S_{3,3}^5})-tE
=
\left(\frac34-t\right)E
+\frac14\widetilde L
+\frac14\widetilde\Gamma_{3,3}.
\]
Its Zariski decomposition is
\[
P(t)=
\begin{dcases}
    \left(\frac34-t\right)
    \left(E+\frac13\widetilde L\right)
    +\frac14\widetilde\Gamma_{3,3},
    & t\in\left[0,\frac12\right],\\
    \left(\frac34-t\right)
    \left(E+\frac13\widetilde L
    +\widetilde\Gamma_{3,3}\right),
    & t\in\left[\frac12,\frac34\right],
\end{dcases}
\]
and
\[
N(t)=
\begin{dcases}
    \frac13t\widetilde L,
    & t\in\left[0,\frac12\right],\\
    \frac13t\widetilde L
    +\left(t-\frac12\right)\widetilde\Gamma_{3,3},
    & t\in\left[\frac12,\frac34\right].
\end{dcases}
\]
Thus, we obtain \(S_{S_{3,3}^5}(E)=\frac5{12}\). Since \(A_{S_{3,3}^5}(E)=\frac12\), we have
\[
    \frac{A_{S_{3,3}^5}(E)}{S_{S_{3,3}^5}(E)}=\frac65.
\]
Moreover,
\[
P(t)\cdot E=
\begin{dcases}
    \frac83t,
    & t\in\left[0,\frac12\right],\\
    \frac{16}{3}\left(\frac34-t\right),
    & t\in\left[\frac12,\frac34\right].
\end{dcases}
\]

For \(z\in E\setminus\left(\widetilde L\cup\Supp(\widetilde\Gamma_{3,3})\right)\),
we have \(S(W_{\bullet,\bullet}^{E};z)=\frac49\) and \(A_{E,\Phi_E}(z)=1\).
For \(z\in E\cap\Supp(\widetilde\Gamma_{3,3})\),
we have \(S(W_{\bullet,\bullet}^{E};z)\leq\frac23\) and \(A_{E,\Phi_E}(z)=1\).
It follows from Lemma \ref{lem:AZ-crepant-model} that
\[
\inf_{z\in E\setminus\{q\}}
\delta_z(S_2^{(2)},B_{S_2^{(2)}};\psi_{3,3}^{\ast}(-K_{S_{3,3}^5}))
\geq
\frac65,
\]
where \(q\coloneqq E\cap\widetilde L\).

At \(q\), we have \(A_{E,\Phi_E}(q)=\frac12\) and \(S(W_{\bullet,\bullet}^{E};q)=\frac7{12}\).
The resulting ratio is \(\frac67\). Hence, this point must be treated on
a higher model.

Let \(\eta\colon Y\to S_2^{(2)}\) be the blow-up of \(q\). Denote its exceptional curve by \(F\), and denote the strict transforms of \(E\) and \(\widetilde L\) by \(E_1\) and \(L_1\), respectively. Put \(\mu\coloneqq\psi_{3,3}\circ\eta\) and \(M_Y\coloneqq\mu^{\ast}(-K_{S_{3,3}^5})\).
Then \(K_Y+B_Y=\mu^{\ast}K_{S_{3,3}^5}\), where \(B_Y=\frac12E_1+\frac12L_1\).

The linear system
\(|\mathcal O_{\P(1,1,3)}(4)|\) has projective dimension \(6\).
Passing through the five general points imposes five linear conditions,
while requiring that the strict transform on \(S^{(2)}_2\) pass through
\(q\) imposes one further linear condition. These conditions are
independent for general points. Let \(\Theta\in|\mathcal O_{\P(1,1,3)}(4)|\) be the resulting curve, let \(T_0\) be its strict transform on
\(S_2^{(2)}\), and let \(T\) be its strict transform on \(Y\). Generality implies that \(T_0\) is smooth at
\(q\) and meets \(E\) and \(\widetilde L\) transversely there.

We have the following intersection numbers
\[
    E_1^2=L_1^2=-4,
    \qquad
    F^2=T^2=-1,
    \qquad
    E_1\cdot F=L_1\cdot F=T\cdot F=1.
\]
The curves \(E_1\), \(L_1\), and \(T\) are pairwise disjoint. Moreover,
\((-K_{S_{3,3}^5})^2=1\), and \(\mu^\ast(-K_{S_{3,3}^5})\equiv 2F+\frac12E_1+\frac12L_1+T\). In particular, \(\mu^\ast(-K_{S_{3,3}^5})\cdot T=1\).

The Zariski decomposition \(\mu^\ast(-K_{S_{3,3}^5})-tF=P(t)+N(t)\) is given by
\[
P(t)=
\begin{dcases}
    (1-t)\mu^\ast(-K_{S_{3,3}^5})+tP_1,
    & t\in[0,1],\\
    (2-t)P_1,
    & t\in[1,2],
\end{dcases}
\qquad
N(t)=
\begin{dcases}
    \dfrac{t}{4}(E_1+L_1),
    & t\in[0,1],\\
    \dfrac{t}{4}(E_1+L_1)+(t-1)T,
    & t\in[1,2],
\end{dcases}
\]
where \(P_1\coloneqq \mu^\ast(-K_{S_{3,3}^5})-F-\frac14(E_1+L_1)\equiv \frac12(\mu^\ast(-K_{S_{3,3}^5})+T)\). Moreover, we have \(\tau(F)=2\).

Thus, we obtain
\[
P(t)^2=
\begin{dcases}
    1-\dfrac12t^2,
    & t\in[0,1],\\
    \dfrac12(2-t)^2,
    & t\in[1,2],
\end{dcases}
\quad
P(t)\cdot F=
\begin{dcases}
    \dfrac12t,
    & t\in[0,1],\\
    \dfrac12(2-t),
    & t\in[1,2].
\end{dcases}
\]
Thus, we obtain \(S_{S_{3,3}^5}(F)=1\). Since \(F\) is obtained by blowing up the transverse intersection of
\(E\) and \(\widetilde L\), we have \(A_{S_{3,3}^5}(F)=A_{S_{3,3}^5}(E)+A_{S_{3,3}^5}(\widetilde L)=\frac12+\frac12=1\).
Hence, 
\[
\frac{A_{S_{3,3}^5}(F)}{S_{S_{3,3}^5}(F)}=1.
\]

Let \(q_E\coloneqq F\cap E_1\), \(q_L\coloneqq F\cap L_1\), \(q_T\coloneqq F\cap T\).
Using \eqref{function h(t)}, we obtain
\[
S(W_{\bullet,\bullet}^F;q)=
\begin{dcases}
    \dfrac16,
    & q\in F\setminus\{q_E,q_L,q_T\},\\
    \dfrac5{12},
    & q\in\{q_E,q_L\},\\
    \dfrac13,
    & q=q_T.
\end{dcases}
\]
The different satisfies
\[
A_{F,\Phi_F}(q)=
\begin{dcases}
    1,
    & q\in F\setminus\{q_E,q_L\},\\
    \dfrac12,
    & q\in\{q_E,q_L\}.
\end{dcases}
\]
In particular,
\[
\frac{A_{F,\Phi_F}(q)}
{S(W_{\bullet,\bullet}^F;q)}
\geq
\min\left\{
    6,\frac65,3
\right\}
=
\frac65>1.
\]

Applying Lemma \ref{lem:AZ-crepant-model} to the crepant pair
\((Y,B_Y)\) and the divisor \(F\), the preceding
calculations give
\[
\inf_{z\in F}
\delta_z(Y,B_Y;M_Y)
\geq
\min\left\{
    \frac{A_{S_{3,3}^5}(F)}{S_{S_{3,3}^5}(F)},
    \inf_{z\in F}
    \frac{A_{F,\Phi_F}(z)}
    {S(W_{\bullet,\bullet}^{F};z)}
\right\}
\geq
\min\left\{
    1,
    \frac65
\right\}
=
1.
\]
Since \(\eta^{-1}(q)=F\), the first assertion of Lemma \ref{lem:AZ-crepant-model} gives \(\delta_q(S_2^{(2)},B_{S_2^{(2)}};\psi_{3,3}^{\ast}(-K_{S_{3,3}^5}))\geq1\).

We have therefore proved \(\delta_z(S_2^{(2)},B_{S_2^{(2)}};\psi_{3,3}^{\ast}(-K_{S_{3,3}^5}))\geq1\) for every \(z\in\psi_{3,3}^{-1}(p)=E\cup\widetilde L\). Another application of Lemma \ref{lem:AZ-crepant-model} gives \(\delta_p(S_{3,3}^5)\geq1\).
On the other hand, the divisor \(F\) gives \(\delta_p(S_{3,3}^5)\leq \frac{A_{S_{3,3}^5}(F)}{S_{S_{3,3}^5}(F)}=1\). Consequently, \(\delta_p(S_{3,3}^5)=1\).
\end{proof}

Now, we show that \(\delta_p(S_{3,3}^5)>1\) for smooth points \(p\). 

\begin{lemma}\label{lem;335 smooth}
Let \(p\in S_{3,3}^5\) be a smooth point. Then we have \(\alpha_p(S_{3,3}^5)\geq\frac{7}{10}\). In particular, \(\delta_p(S_{3,3}^5)>1\).
\end{lemma}

\begin{proof}
Apply Lemma \ref{lem:smooth-alpha-reduction} with
\((n,m,k)=(3,3,5)\). The target of the elementary transformation is
\(S_{2,3}^4\).

We identify this target with \(S_{3,2}^5\). The symmetry of the
construction recalled before Figure \ref{fig:two negative} gives \(S_{3,2}^3\cong S_{2,3}^2\).
Let \(\theta_0\colon S_{2,3}^2\to S_{3,2}^3\)
be the isomorphism obtained from the symmetry of the construction.
Let \(u_1,u_2\in S_{2,3}^2\) be the two chosen general smooth points
whose blow-ups define \(S_{2,3}^4\). Put \(v_i\coloneqq\theta_0(u_i)\) for \(i\in\{1,2\}\).
The points \(v_1,v_2\) are smooth and general on \(S_{3,2}^3\).
Blowing up the corresponding pairs of points gives an isomorphism \(\theta\colon S_{2,3}^4\to S_{3,2}^5\).
Consequently, for every smooth point \(q\in S_{2,3}^4\), we have
\[
\alpha_q(S_{2,3}^4)
=
\alpha_{\theta(q)}(S_{3,2}^5)
\geq
\frac7{10}
\]
by Lemma \ref{lem;325 smooth}. Lemma
\ref{lem:smooth-alpha-reduction} now gives
\(\alpha_p(S_{3,3}^5)\geq\frac{7}{10}\).

By Theorem \ref{thm:alpha delta}, we have
\[
    \delta_p(S_{3,3}^5)
    \geq
    \frac32\alpha_p(S_{3,3}^5)
    \geq
    \frac{21}{20}
    >
    1.\qedhere
\]
\end{proof}

\begin{theorem}\label{thm:335-semistable}
The surface \(S_{3,3}^5\) is strictly K-semistable.
\end{theorem}

\begin{proof}
Lemmas \ref{lem;335 sing} and \ref{lem;335 smooth} give \(\delta(S_{3,3}^5)=1\). Hence, \(S_{3,3}^5\) is K-semistable by Theorem \ref{thm:delta}, and it
is not K-stable by Theorem \ref{thm:K-stable-uniform}. Its automorphism
group is finite by Lemma \ref{lem:335-527-aut-fin}. Consequently,
Theorem \ref{thm:finite aut} shows that it is not K-polystable.
Therefore, \(S_{3,3}^5\) is strictly K-semistable.
\end{proof}

\subsubsection{Case: \(S_{4,3}^6\)}

\begin{lemma}\label{lem;436 sing}
Let \(p\in S_{4,3}^6\) be the singular point of type
\(\frac1{11}(1,4)\). Then we have \(\delta_p(S_{4,3}^6) \leq \frac{27}{29} <1\).
\end{lemma}

\begin{proof}
Let \(q\coloneqq E\cap\widetilde L\), and let
\(\eta\colon Y\to S^{(2)}_2\) be the blow-up of \(q\). Denote the
exceptional curve by \(F\), and denote the strict transforms of
\(E\) and \(\widetilde L\) by \(E_1\) and \(L_1\), respectively. Put
\(\mu\coloneqq\psi_{4,3}\circ\eta\).

The linear system
\(|\mathcal O_{\P(1,1,4)}(5)|\) has projective dimension \(7\).
Passing through the six general points imposes six linear conditions,
while requiring that the strict transform on \(S^{(2)}_2\) pass through
\(q\) imposes one further linear condition. These conditions are
independent for general points. Let \(\Theta\in|\mathcal O_{\P(1,1,4)}(5)|\) be the resulting curve, let \(T_0\) be its strict transform on
\(S_2^{(2)}\), and let \(T\) be its strict transform on \(Y\). Generality implies that \(T_0\) is smooth at
\(q\) and meets \(E\) and \(\widetilde L\) transversely there.

We have the following intersection numbers
\[
    E_1^2=-5,
    \qquad
    L_1^2=-4,
    \qquad
    F^2=T^2=-1,
    \qquad
    E_1\cdot F=L_1\cdot F=T\cdot F=1.
\]
The curves \(E_1\), \(L_1\), and \(T\) are pairwise disjoint. Moreover,
\((-K_{S_{4,3}^6})^2=\frac9{11}\), and 
\[
    \mu^\ast(-K_{S_{4,3}^6})
    \equiv
    \frac{20}{11}F+\frac4{11}E_1+\frac5{11}L_1+T.
\]
In particular, we have \(\mu^\ast(-K_{S_{4,3}^6})\cdot T=\frac9{11}\).

The Zariski decomposition of \(\mu^\ast(-K_{S_{4,3}^6})-tF\) is given by
\[
P(t)=
\begin{dcases}
    \left(1-\dfrac{11t}{9}\right)\mu^\ast(-K_{S_{4,3}^6})
    +\dfrac{11t}{9}P_\ast,
    & t\in\left[0,\frac9{11}\right],\\
    \left(\dfrac{20}{11}-t\right)P_\ast,
    & t\in\left[\frac9{11},\frac{20}{11}\right].
\end{dcases}
\]
and
\[
N(t)=
\begin{dcases}
    \dfrac{t}{5}E_1+\dfrac{t}{4}L_1,
    & t\in\left[0,\frac9{11}\right],\\
    \dfrac{t}{5}E_1+\dfrac{t}{4}L_1
    +\left(t-\dfrac9{11}\right)T,
    & t\in\left[\frac9{11},\frac{20}{11}\right],
\end{dcases}
\]
where \(P_\ast \coloneqq \frac{11}{20}\mu^\ast(-K_{S_{4,3}^6})+\frac9{20}T\).
The divisor \(P_\ast\) is nef. Indeed, its intersections with
\(E_1\), \(L_1\), and \(T\) are zero, while
\(P_\ast\cdot F=\frac9{20}\). For every other irreducible curve \(C\),
we have \(\mu^\ast(-K_{S_{4,3}^6})\cdot C>0\) and \(T\cdot C\geq0\). 
Moreover, \(\tau(F)=\frac{20}{11}\).

We have
\[
P(t)^2=
\begin{dcases}
    \dfrac9{11}-\dfrac{11}{20}t^2,
    & t\in\left[0,\frac9{11}\right],\\
    \dfrac9{20}\left(\dfrac{20}{11}-t\right)^2,
    & t\in\left[\frac9{11},\frac{20}{11}\right].
\end{dcases}
\]
Thus, we obtain \(S_{S_{4,3}^6}(F)=\frac{29}{33}\). Since \(F\) is obtained by blowing up the transverse intersection of
\(E\) and \(\widetilde L\), we have \(A_{S_{4,3}^6}(F) = A_{S_{4,3}^6}(E)+A_{S_{4,3}^6}(\widetilde L) = \frac4{11}+\frac5{11} = \frac9{11}\).
Hence, we obtain 
\[
\frac{A_{S_{4,3}^6}(F)}{S_{S_{4,3}^6}(F)} = \frac{27}{29} < 1.
\]
Since \(F\) is a prime divisor over \(S_{4,3}^6\), it follows that \(\delta_p(S_{4,3}^6) \leq \frac{27}{29} <1\).
Equivalently, 
\[
\beta_{S_{4,3}^6}(F) = \frac9{11}-\frac{29}{33} = -\frac2{33} <0.\qedhere
\]
\end{proof}

\begin{lemma}\label{lem;436 smooth}
Let \(p\in S_{4,3}^6\) be a smooth point. Then we have \(\alpha_p(S_{4,3}^6)\geq\frac{7}{10}\). In particular, \(\delta_p(S_{4,3}^6)>1\).
\end{lemma}

\begin{proof}
Apply Lemma \ref{lem:smooth-alpha-reduction} with
\((n,m,k)=(4,3,6)\). The target of the elementary transformation is
\(S_{3,3}^5\). By Lemma \ref{lem;335 smooth}, we obtain \(\alpha_q(S_{3,3}^5)\geq\frac{7}{10}\) for every smooth point \(q\in S_{3,3}^5\). Hence,
\(\alpha_p(S_{4,3}^6)\geq\frac{7}{10}\).

By Theorem \ref{thm:alpha delta}, we have
\[
    \delta_p(S_{4,3}^6)
    \geq
    \frac32\alpha_p(S_{4,3}^6)
    \geq
    \frac{21}{20}
    >
    1.\qedhere
\]
\end{proof}

\begin{theorem}\label{thm:436-unstable}
The surface \(S_{4,3}^6\) is K-unstable.
\end{theorem}

\begin{proof}
By Lemma \ref{lem;436 sing}, we have \(\delta(S_{4,3}^6) \leq \frac{27}{29} <1\). Therefore, \(S_{4,3}^6\) is K-unstable by Theorem \ref{thm:delta}.
\end{proof}

Finally, we prove the main result of this paper.
\begin{proof}[Proof of Theorem \ref{main theorem}]
The assertion follows from Remark \ref{rem:22} and
Theorems \ref{thm:solution set},
\ref{thm:325 stable},
\ref{thm:326 stable},
\ref{thm:426 stable},
\ref{thm:427 stable},
\ref{thm:335-semistable},
\ref{thm:436-unstable}, and
\ref{thm:527-semistable}.
\end{proof}

\bibliographystyle{habbvr}
\bibliography{biblio}

\end{document}